\documentclass[a4paper]{amsart}
\usepackage{graphicx}
\usepackage{amssymb}
\usepackage{amsfonts}
\usepackage{hyperref}
\usepackage{mathrsfs}

\usepackage[dvipsnames]{xcolor}
\swapnumbers
\newtheorem{theorem}{Theorem}[section]
\newtheorem{lemma}[theorem]{Lemma}
\newtheorem{corollary}[theorem]{Corollary}
\newtheorem{proposition}[theorem]{Proposition}
\theoremstyle{definition}
\newtheorem{definition}[theorem]{Definition}
\newtheorem{assumption}[theorem]{Assumption}
\newtheorem{remark}[theorem]{Remark}
\newtheorem{example}[theorem]{Example}

\numberwithin{equation}{section}
\theoremstyle{plain}

\numberwithin{equation}{section} 
\numberwithin{figure}{section} 
\theoremstyle{plain}
\theoremstyle{plain}
\theoremstyle{remark}
\newtheorem*{acknowledgement*}{Acknowledgement}

\newcommand{\cB}{{\mathcal B}}

\newcommand{\cD}{{\mathcal D}}
\newcommand{\cE}{{\mathcal E}}
\newcommand{\cF}{{\mathcal F}}

\newcommand{\cT}{{\mathcal T}}

\newcommand{\cX}{{\mathcal X}}

\newcommand{\Om}{{\Omega}}
\newcommand{\om}{{\omega}}

\newcommand{\bbE}{{\mathbb E}}

\newcommand{\bbP}{{\mathbb P}}
\newcommand{\bbR}{{\mathbb R}}

\newcommand{\bbZ}{{\mathbb Z}}

\AtBeginDocument{\bibliographystyle{amsplain}}
\begin{document}
\title[]{Local Limit Theorems and Edgeworth Expansions for some classes of inhomogeneous Markov chains}  
 \vskip 0.1cm
 \author{Yeor Hafouta$^1$, Lubashan Pathirana$^2$, Brendan Williams$^3$
}
\thanks{1. Department of Mathematics, Ben-Gurion University of the Negev, Israel. \,\,\,Email: yeorh@bgu.ac.il}
	\qquad
	\thanks{2. Department of Mathematical Sciences, University of Copenhagen, Denmark.\,\,\,
	Email: lpk@math.ku.dk}
    \qquad
    \thanks{3. Department of Mathematics, University of Florida, USA. \,\,\,
	Email: brendan.williams@ufl.edu }


\thanks{ }
\dedicatory{  }
 \date{\today}

\begin{abstract}
    In this note, we prove local central limit theorems and expansions in the central limit theorem for inhomogeneous Markov chains satisfying a one-step lower ellipticity condition. As for the local limit theorem, in the general inhomogeneous case, compared with \cite{DS}, we are able to prove results without two-sided ellipticity, although only in the one-step case. Compared with \cite{PeligradLLT1,PeligradLLT2}, we are able to prove local limit theorems without a lower $\psi$-mixing condition. Concerning Edgeworth expansions, compared with \cite{DHPTRF1}, we are able to get results without two-sided ellipticity (under some sequential Diophantine conditions). Our conditions are then verified for Markov chains in a random dynamical environment with a one-sided Doeblin condition under very mild assumptions. For Markov chains in a random environment, we also consider the more general case of $\phi$-mixing chains and then, under some mild periodicity assumptions on the environment, we prove local central limit theorems and first-order Edgeworth expansions  (without any ellipticity assumptions), matching the results for deterministic homogeneous chains. When the chain is driven by a topologically mixing subshift of finite type, we obtain high-order expansions for $\phi$-mixing chains under some Diophantine conditions of a deterministic system generated by a periodic point.
    We then apply our results to bounded additive functionals of posterior states in disordered quantum trajectories and establish local central limit theorems and first-order Edgeworth expansions.
\end{abstract}

\maketitle
\markboth{Y. Hafouta, L. Pathirana, B. Williams }{ } 
\renewcommand{\theequation}{\arabic{section}.\arabic{equation}}
\pagenumbering{arabic}

\section{Introduction}
One of the most celebrated results in probability theory is the central limit theorem (CLT). In its basic form, it states that for an independent square-integrable zero-mean sequence of random variables $Y_1,Y_2,...$ satisfying the Lindeberg condition, the sequence $S_n/\|S_n\|_{L^2}$, $S_n=Y_1+\ldots +Y_n$ weakly converges towards the standard normal law.   In fact, the first appearance of the CLT in literature is in the so-called De Moivre-Laplace theorem, whose modern version is the local central limit theorem (LLT). The LLT 
concerns the asymptotics of probabilities of the form $\mathbb P(S_n\in I)$ for an interval $I$. 

Once the CLT is achieved, it is desirable to see when the optimal uniform rates $O(\|S_n\|_{L^2}^{-1})$ can be obtained.
Other finer versions of the CLT include
Edgeworth expansions, which provide a more accurate approximation (better than $O(\|S_n\|_{L^2}^{-1})$)  of the distribution function of $S_n/\|S_n\|_{L^2}$. In the past century, the LLT and Edgeworth expansions have been studied extensively for stationary sequences $X_j$ generated by either a sufficiently fast mixing homogeneous Markov chain or a chaotic dynamical system; see, for instance, \cite{FL,HH} and references therein for general approaches.

Stationarity is a very strong assumption. Many real-life models are non-stationary due to interaction with the outside world (e.g., external forces). 
In this paper, we focus on partial sums of functionals $Y_j=f_j(X_j)$ of inhomogeneous Markov chains. The CLT in that setting was first proved by Dobrushin \cite{Dob,Dob2} under a certain contraction condition. We also refer to \cite{SV} for a modern approach to Dobrushin's CLT and to \cite{PeligradPTRF} for a generalization to certain classes of $\rho$-mixing Markov chains.

Recently, the LLT and Edgeworth expansions for inhomogeneous Markov chains were  studied in \cite{DS, DHPTRF1,PeligradLLT1,PeligradLLT2}. See also \cite{HW2} for Berry-Esseen theorems under $\phi$-mixing. The results in \cite{DHPTRF1, DS} assumed a uniform two-sided ellipticity assumption, while the results in \cite{PeligradLLT1,PeligradLLT2} assumed the lower $\psi$-mixing condition, and they did not include Edgeworth expansions. In this note, we will obtain an LLT and Edgeworth expansions under a lower one-step ellipticity condition (a lower Doeblin condition). Then we apply these results to Markov chains in random dynamical environments. In the latter random case,  
we also present a second approach and, under mild assumptions on the environment, we obtain results under $\phi$-mixing conditions, which include chains satisfying a lower Doeblin condition of an arbitrary order or a uniform Dobrushin-type contraction. We will also discuss certain applications to disordered quantum trajectories.

Our approach in the inhomogeneous (not random) case relies on a simple general upper bound on the characteristic function of the underlying partial sum which holds for any argument $t\in\mathbb R$ (see Lemma \ref{Key}). This estimate is in the spirit of  \cite[Lemma 11]{PeligradLLT2}, but with a different norm (and a shorter proof). It also leads to a shorter proof of a version of \cite[Proposition 10]{PeligradLLT2}. 
Then, we use the results in \cite{HW2} and ideas in \cite{DS} to get the desired decay of the characteristic functions needed for the LLT and the expansions. For one-step lower elliptic Markov chains in random environment, our conditions for the LLT hold under a random irreducibility condition, and the conditions for the Edgeworth expansions hold under a random Diophantine condition. In the general $\phi$-mixing case in random environment, we use ideas in \cite[Chapter 7]{HK} to get certain bounds on the characteristic function on long intervals, which yield the desired results.

\section{Main results for inhomogeneous one-step lower elliptic Markov chains}
Let $(X_j)_{j=1}^\infty$ be a Markov chain with state spaces $(\cX_j,\cF_j)$ and transition probabilities $R_j(x,\Gamma)=\mathbb P(X_{j+1}\in\Gamma|X_j=x)$, $x\in\cX_j, \Gamma\in\cF_{j+1}$. We will work under the following assumption.
\begin{assumption}[One-step lower ellipticity]\label{A1}
There exist probability measures $m_j$ on $(\cX_{j},\cF_j)$ and a constant $\gamma\in(0,1)$ such that for all $j\geq1$, $x\in\cX_j$ and $\Gamma\in\cF_{j+1}$, we have 
$$
R_j(x,\Gamma)\geq \gamma m_{j+1}(\Gamma).
$$
\end{assumption}
 Notice that under Assumption \ref{A1} we have 
\begin{equation}\label{Dob}
\sup_j\sup_{x,y\in\cX_j}\sup_{\Gamma\in\cF_{j+1}}\left|R_j(x,\Gamma)-R_j(y,\Gamma)\right|\leq 1-\gamma.    
\end{equation}
That is, we have a one-step Dobrushin contraction.

Next, let $f_j:\cX_j\to\mathbb R$ be uniformly bounded measurable functions, and let 
$$
S_n=\sum_{j=1}^n f_j(X_j).
$$
Denote $V_n=\text{Var}(S_n)$. Then by \eqref{Dob} and
\cite[Proposition 13]{PeligradPTRF} there are constants $C_1,C_2>0$ which depend only on $\gamma$ such that for all $n\in\mathbb N$,
\begin{equation}
\label{eq:variance-comparison}
    C_1\sum_{j=1}^n\text{Var}(f_j(X_j))\leq V_n\leq C_2\sum_{j=1}^n\text{Var}(f_j(X_j)).
\end{equation}
 When $\sum_{j=1}^\infty\text{Var}(f_j(X_j))=\infty$
under Assumption \ref{A1} we showed in \cite{HW2} that 
$V_n^{-1/2}(S_n-\mathbb E[S_n])$ weakly converges to the standard normal distribution with optimal uniform rates. That is, a Berry--Esseen theorem holds.
In this paper, we will complement this by proving a local limit theorem and Edgeworth expansions.
\begin{definition}
\label{dfn:LCLT}
We say that $S_n$ obeys the non-lattice local central limit theorem if $V_n\to\infty$ and  for every continuous function $G:\mathbb R\to\mathbb R$ with compact support or an indicator of a bounded interval, we have 
$$
\sup_{u\in\mathbb R}\left|\sqrt{2\pi V_n}\mathbb E[G(S_n-u)]-\left(\int_{\mathbb R}G(x)dx\right)e^{-\frac{(u-\mathbb E[S_n])^2}{2V_n}}\right|=o(1).
$$
\end{definition}

\begin{definition}
\label{dfn:1_edgeworth}
We say that $S_n$ obeys the first-order Edgeworth expansions\footnote{Yeor I changed the polynomial} if $V_n\to\infty$   and 
$$
\sup_{t\in\mathbb R}\left|\mathbb P(S_n-\mathbb E[S_n]\leq t\sqrt{V_n})-\Phi(t)-V_n^{-1/2}\frac{1}{6}(1-t^2)\gamma_{3,n}\varphi(t)\right|=o(V_n^{-1/2})
$$
where $\varphi(t)=\frac1{\sqrt{2\pi}}e^{-\frac 12 t^2}$, $\Phi(t)=\int_{-\infty}^t \varphi(x)dx$ and $\gamma_{3,n}=V_n^{-1}\mathbb E[(S_n-\mathbb E[S_n])^3]$. 
\end{definition}

We note that in the setting of all the results in this paper, the term $\gamma_{3,n}$ is uniformly bounded in $n$, see \cite[Proposition 3.2]{HW2}.

Our first result is.
\begin{theorem}\label{LLT0}
Make Assumption \ref{A1}. Suppose also $\sum_{j=1}^\infty\text{Var}(f_j(X_j))=\infty$. 
For $j\geq0$ and $t\in\mathbb R$ let
$$
d_j(t)=1-\left|\int e^{itf_j(y)}dm_j(y)\right|.
$$
Then the non-lattice local central limit theorem and the first-order Edgeworth expansions hold if for every compact set $K\subset \mathbb R\setminus\{0\}$ we have
\begin{equation}\label{llt}
\lim_{n\to\infty}V_n^{1/2}\int_K e^{-\frac12\gamma^2\sum_{j=1}^n d_j(t)}dt=0.    
\end{equation}
\end{theorem}
Theorem \ref{LLT0} is in the spirit of the LLT's in \cite{PeligradLLT1,PeligradLLT2}, but with $m_j$ replacing the law of $X_j$.
Condition \eqref{llt} in Theorem \ref{LLT0} means that most of $f_j$ are not close to a lattice-valued random variable with jumps $2\pi/t$ (with respect to $m_j$).

Under certain assumptions on the relation between $m_j$ and the distribution of $X_j$, we are able to prove a local limit theorem without the decay rates in \eqref{llt}.
\begin{theorem}\label{LLT1}
Make Assumption \ref{A1}. Suppose also $\sum_{j=1}^\infty\text{Var}(f_j(X_j))=\infty$. Denote $v_n=\sum_{j=1}^n\text{Var}_{m_j}(f_j)$ and assume that $v_n\geq c_1\sum_{j=1}^n\text{Var}(f_j(X_j))$ for some constant $c_1>0$ and all $n$ large enough.
Assume also that for all nonzero real $t$, we have 
\begin{equation}\label{div}
\sum_{j=1}^{\infty}\left(1-\left|\int e^{itf_j(y)}dm_j(y)\right|\right)=\infty.    
\end{equation}
Then the non-lattice local central limit theorem and the first-order Edgeworth expansions hold.
\end{theorem}
Theorem \ref{LLT1} is new already under the assumption that $m_j$ is (equivalent to) the law of $X_j$ like in  \cite{PeligradLLT1,PeligradLLT2}. This case corresponds to one-step lower $\psi$-mixing. The condition $v_n\geq c_1\sum_{j=1}^n\text{Var}(f_j(X_j))$ clearly holds when $\text{Var}_{m_j}(f_j)\geq c$ for some constant $c>0$ and all $j$ large enough. In Theorem \ref{LLT rand0} for Markov chains in a random environment, it is clear that both $v_n$ and $V_n$ grow linearly fast in $n$, and so this condition holds.
However, in general $\sum_{j=1}^n\text{Var}(f_j(X_j))$ can grow much slower than $n$ and then decay of $\text{Var}_{m_j}(f_j)$ is allowed. For instance, if $\text{Var}(f_j(X_j))=O(j^{-\beta})$ for $0<\beta\leq 1$ then we only require $v_n\geq cn^{1-\beta}$ for $\beta<1$ and $v_n\geq c\ln n$ for $\beta=1$, which are quite mild growth conditions.

To see why the condition \eqref{div} is necessary, we have the following result.
\begin{theorem}\label{LLT Red}
Make Assumption \ref{A1}.
Suppose that for some $t\not=0$
$$
\sum_{j=1}^{\infty}\left(1-\left|\int e^{itf_j(y)}dm_j(y)\right|\right)<\infty.
$$  
Then for all $j$ we have 
\begin{equation}\label{Red0}
f_j=a_j+(2\pi/t)Z_j+U_j,\,\,\,m_j-\text{a.s.}
\end{equation}
where the $a_j\in\mathbb R$, $Z_j$ is an integer-valued function and $\sum_j\int_{\cX_j}U_j^2\,dm_j<\infty$. Moreover,  on the product space $\prod_{j}(\cX_j,\cF_j ,m_j)$ we have that 
\begin{equation}\label{Red}
f_j=a_j+(2\pi/t)Z_j+U_j    
\end{equation}
where the sum $\sum_j U_j$ converges almost surely and in $L^2$.
\end{theorem}

In the above circumstances, it is not expected that the non-lattice local CLT would hold, even when $m_j$ is the law of $X_j$. Indeed, in that case $f_j$ can be reduced to a lattice-valued sequence. Thus, condition \eqref{div} in Theorem \ref{LLT1} is necessary.
\begin{remark}\label{r1}
The only place where the (uniform) boundedness of $f_j$ is used in the proof of the local limit theorem in  Theorem \ref{LLT1} is to ensure that the central limit theorem holds (see Remark \ref{bb}).
This can be obtained under some growth conditions on $\|f_j\|_\infty$; see \cite{PeligradPTRF}. Assuming the validity of the central limit theorem, our proof proceeds similarly if, for instance, we assume that $\sup_j\|f_j\|_{L^3(m_j)}<\infty$.
However, we decided to state the main results under uniform boundedness for the sake of simplicity.
\end{remark}

Next, we discuss high-order Edgeworth expansions.
\begin{definition}
Let $r\in \mathbb N$, $r\geq2$.
We say that $S_n$ obeys Edgeworth expansions of order $r$ if there are polynomials $P_{j,n}$ with degrees depending only on $j$ and with uniformly bounded in $n$ coefficients such that 
$$
\sup_{t\in\mathbb R}\left|\mathbb P(S_n-\mathbb E[S_n]\leq t\sqrt{V_n})-\Phi(t)-\sum_{j=1}^rV_n^{-j/2}P_{j,n}(t)\varphi(t)\right|=o(V_n^{-r/2}).
$$
\end{definition}
\begin{theorem}\label{Edge1}
Make Assumption \ref{A1}. Let $r\geq2$ and suppose that for all $A,B>0$, we have
\begin{equation}\label{EdgeDec}
\lim_{n\to\infty}V_n^{r/2}\int_{A\leq |t|\leq BV_n^{(r-1)/2}}|t|^{-1}e^{-\frac12\gamma^2\sum_{j=1}^n d_j(t)}dt=0,    
\end{equation}
where $d_j(t)=1-\left|\int e^{itf_j(y)}dm_j(y)\right|$.
 Then $S_n$ obeys Edgeworth expansions of order $r$.
\end{theorem}

The next result shows that 
condition \eqref{EdgeDec} holds under uniform Diophantine conditions.
\begin{proposition}\label{EdgeProp}
Condition \eqref{EdgeDec} holds if 
there are $c,\delta_0,\ell>0$  such that $\frac12(r-1)\ell<1$  and for all $t\in\mathbb R$ with $|t|\geq\delta_0$, we have
$$
\sup_j\left|\int e^{itf_j(y)}dm_j(y)\right|\leq 1-c|t|^{-\ell}.
$$
\end{proposition}

\begin{remark}
Our proof shows that
condition \eqref{EdgeDec} also holds when there are  $c_j>0$ such that
$$
\left|\int e^{itf_j(y)}dm_j(y)\right|\leq 1-c_j|t|^{-\ell}
$$
and  $\sum_{j=1}^n c_j\geq cn^{\kappa}$ for some $\kappa>0$ such that $\ell(r-1)/2<\kappa$. This allows us to obtain expansions under non-uniform in $j$ Diophantine conditions.
\end{remark}

\section{Main results for Markov chains in random environment}
Let $(\Omega,\cF,\bbP)$ be a  probability space. Let $\theta:\Omega\to\Omega$ be an ergodic invertible probability-preserving transformation. Let $\cX$ be a measurable space and let us take a set $\cE\subset \Omega\times \cX$ measurable with respect to the product $\sigma$-algebra. Let $\cE_\omega=\{x\in\cX: (\omega,x)\in\cE\}$. For $\omega\in\Omega$ let $(X_{\omega,k})_{k\in\mathbb Z}$ be a Markov chain with state spaces $\cE_{\theta^k\omega}$ and transition probability $R_{\theta^k\omega}(x,dy)$ from time $k$ to time $k+1$ depending only on  $\theta^k\omega$. 

Next, let $f:\Omega\times\cX\to\mathbb R$ be a bounded measurable function. Denote $f_\omega(x)=f(\omega,x)$ and set 
$$
S_n^\omega=\sum_{j=1}^{n}f_{\theta^j\omega}(X_{\omega,j}).
$$

Our first result is.
\begin{theorem}\label{LLT rand0}
Suppose that there are probability measures $m_\omega$ on $\mathcal E_\omega$ measurable in $\omega$ and $\gamma\in(0,1)$ such that for $\mathbb P$-a.a. $\omega$ for all $x\in\cE_{\omega}$, we have 
$$
R_{\om}(x,\cdot)\geq \gamma m_{\theta\omega}(\cdot).
$$
Suppose also that $\text{Var}(S_n^\omega)\to\infty$ almost surely and that $\omega\to \text{Var}_{m_\omega}(f_\omega)$ is not the $0$-random variable
\vskip0.1cm
(i) Assume also that for all $h\not=0$ we have 
$$
\mathbb P(\Gamma_h)<1
$$
where\footnote{Note that $\Gamma_h$ is also the collection of all $\omega\in\Omega$ such that the function $(x,y)\to f_\omega(x)-f_\omega(y)$ takes values in $h\mathbb Z$ for $m_\omega\times m_\omega$-a.a. $(x,y)$.}
$\Gamma_h=\{\omega\in\Omega: |m_\omega(e^{it f_\omega})|=1\}$, for $t=\frac{2\pi}{h}$.
Then $\mathbb P$-a.a. $\omega$, the sequence $S_n^\omega$ obeys the non-lattice local central limit theorem and the first-order Edgeworth expansion. \\
\vskip0.1cm
(ii) Let $r\geq 2$. Assume that with positive probability there exist $c_\omega>0$ and $\delta_\omega>0$ such that whenever $|t|>\delta_\omega$ we have $|m_\omega(e^{itf_\omega})|\leq 1-c_\omega |t|^{-\ell}$ with $\ell>0$ satisfying $(r-1)\ell<2$ (i.e. a random Diophantine condition). Then for $\mathbb P$-a.a. $\omega$, the sequence $S_n^\omega$ obeys the Edgeworth expansions of order $r$.
\end{theorem}

The proof of Theorem \ref{LLT rand0} is based on verifying the conditions of Theorems \ref{LLT1} and \ref{Edge1} and the ergodicity of $\theta$.

\begin{remark}
 Note that by  \cite[Theorem 2.1]{Kifer1998}, there is a number $\sigma\geq0$ such that for $\mathbb P$-a.a. $\omega$ we have 
 $$
 \sigma^2=\lim_{n\to\infty}\frac1n\text{Var}(S_n^\omega).
$$ 
Thus the condition $\text{Var}(S_n^\omega)\to\infty$ means that $\sigma>0$. Note also that by \eqref{eq:variance-comparison} we have 
 $$
C_1\sum_{j=1}^n\text{Var}(f_{\theta^j\omega}(X_{\omega,j}))\leq \text{Var}(S_n^\omega)\leq C_2\sum_{j=1}^n\text{Var}(f_{\theta^j\omega}(X_{\omega,j})).
 $$
Moreover, in \cite{HW2} we showed that there is a family of probability measures $\mu_\omega$ and constants $C>0$ and $\delta\in(0,1)$ such that 
 $R_{\omega}^*\mu_{\omega}=\mu_{\theta\omega}$ (a.s.) and $\mathbb P$-a.s. for all $n\in\mathbb N$,
$$
\sup_{\|g\|_\infty\leq 1}\left\|R_{\omega}\circ\cdots\circ R_{\theta^{j-1}\omega}g-\mu_{\theta^j\omega}(g)\right\|_\infty\leq C\delta^j.
$$
 Therefore, the difference between the law of $X_{\omega,j}$ and $\mu_{\theta^j\omega}$ converges to $0$ exponentially fast in total variation as $j\to\infty$ (uniformly in $\omega$). By using the mean ergodic theorem with the function $c(\omega)=\text{Var}_{\mu_\omega}(f_\omega)$ we
 conclude that 
 $$
C_1\int c(\omega)d\mathbb P(\omega)\leq \sigma^2\leq C_2\int c(\omega)d\mathbb P(\omega).
 $$
 Hence $\sigma=0$ if and only if $\text{Var}_{\mu_\omega}(f_\omega)=0$ a.s. Thus, the condition $\text{Var}(S_n^\omega)\to\infty$ means that $f_\omega$ is not constant $\mu_\omega$-a.s. for $\mathbb P$-a.a. $\omega$. 
\end{remark}

\begin{remark}\label{r2}
By  \cite[Theorem 2.1]{Kifer1998}, when $\sigma>0$ the central limit theorem for $S_n^\omega$ holds when $\|f_\omega\|_{L^p(\mu_\omega)}\in L^\infty(\Omega,\mathcal F,\mathbb P)$ for some $p>2$. The proof of the local central limit theorem in Theorem  \ref{LLT rand0} (i) proceeds the same under this assumption and the assumption that $\|f_\omega\|_{L^3(m_\omega)}\in L^\infty(\Omega,\mathcal F,\mathbb P)$ instead of uniform boundedness of $f$. See Remark \ref{bb}.
\end{remark}

\begin{remark}
\label{rem:nonlattice-criteria}
    In the setting of Theorem~\ref{LLT rand0} above, take $\lambda_\om = (f_\om)_*m_\om$. Then the non-lattice condition required in part~(i) holds if there exists an event $B\subseteq\Omega$ with $\mathbb P(B)>0$ such that, for every $\omega\in B$,  the support of $\lambda_\om$ is infinite or contains three points $a_1, a_2, a_3$ with $a_1\neq a_2$ and 
    \[
        \frac{a_3-a_1}{a_2-a_1}\notin\mathbb Q.
    \]
    These points may depend on $\omega$.
    Indeed, since $f$ is bounded, we have that $\operatorname{supp}(\lambda_\om)$ is also bounded and thus there is some $M>0$ such that $\operatorname{supp}(\lambda_\om)\subseteq [-M,M]$ for every $\om\in\Om$. 
    Now, for any $\om\in B$, if $\lambda_\om$ were supported on a lattice $a+h\mathbb Z$ (for $h>0$), its support would be contained in the finite set $[-M, M]\cap (a+h\mathbb Z)$, excluding the infinite support case. 
    Also, for any three lattice points $b_1, b_2, b_3$ with $b_1\neq b_2$, the ratio $(b_3 - b_1)/(b_2-b_1)$ is rational.
    This means that $\lambda_\omega$ is non-lattice for every $\omega\in B$.
    These criteria imply the stronger property that the same positive $\mathbb P$-probability event excludes every lattice spacing $h>0$.  
\end{remark}

\medskip

Now let us consider more general chains. Let $\phi_\omega(n)$ be the $\phi$-mixing coefficients of the chain $(X_{\omega,k})_{k\in\mathbb Z}$. That is, denote by $\kappa_\omega$ the law of the path $(X_{\omega,k})_{k\in\mathbb Z}$. Then
$$
\phi_{\omega}(n)=\sup_k\sup\left\{|\kappa_\omega(B|A)-\kappa_\omega(B)|:\kappa_\omega(A)>0, A\in\sigma_{\omega,-\infty,k}, B\in\sigma_{\omega,k+n,\infty}\right\},
$$
where for all $-\infty\leq k\leq\ell\leq\infty$, we have $\sigma_{\omega,k,\ell}=\sigma\{X_{\omega,s}: k\leq s\leq \ell, s\in\mathbb R\}$.
\begin{assumption}\label{A2}
There exist $C>0$ and $\delta\in(0,1)$ such for every $\omega\in\Omega$ for all $n$ we have 
$$
\phi_\omega(n)\leq C\delta^n.
$$
\end{assumption}
Two main examples are chains satisfying some finite-step lower Doeblin condition and finite-step Dobrushin-type contraction; see \cite{HW2}. In particular, Assumption \ref{A1} is stronger than Assumption \ref{A2}. 
\begin{remark}\label{rem}
 Assumption \ref{A2} together with \cite[Theorem 2.1]{Kifer1998} ensure that there exists a number $\sigma\geq0$ such that for $\mathbb P$-a.a. $\omega$ we have $\sigma^2=\lim_{n\to\infty}\frac1n\text{Var}(S_n^\omega)$.  
\end{remark}

Next, we work under the following assumption on $(\Omega,\mathcal F,\mathbb P)$ and the maps $\omega\to R_\omega$ and $\omega\to f_\omega$.
\begin{assumption}\label{A3}
\,
\begin{enumerate}
\item[(i)] $\Omega$ is a topological space and $\mathcal F$ is the Borel $\sigma$-algebra.
\vskip0.1cm
\item[(ii)] there exists a $\theta$-periodic point $\omega_0\in\Omega$.
\vskip0.1cm
\item[(iii)] the measure $\mathbb P$ assigns positive mass to all open neighborhoods of the points in
the finite periodic orbit $\{\theta^j\omega_0: j\geq0\}$.
\vskip0.1cm
\item[(iv)] the map $\theta$ is continuous at the finite periodic orbit $\{\theta^j\omega_0: j\geq0\}$.
\vskip0.1cm
\item[(v)] the map $\omega\to f(\omega,\cdot)$ is continuous with respect to the supremum norm at the finite periodic orbit $\{\theta^j\omega_0: j\geq0\}$.
\vskip0.1cm
\item[(vi)] $\mathcal E_\omega$ is constant in $\omega$ around points in the finite periodic orbit $\{\theta^j\omega_0: j\geq0\}$ and the map $\omega\to R_{\omega}$ is continuous with respect to the supremum operator norm at the finite periodic orbit $\{\theta^j\omega_0: j\geq0\}$.
\end{enumerate}
\end{assumption}
\begin{example}
\label{ex:periodic-subshift}
Assumptions \ref{A3} (i)-(iv) hold when $(\Omega,\mathcal F,\mathbb P,\theta)$ is a topological subshift of finite type \cite{Bowen} and $\mathbb P$ is a Gibbs measure. In this case, periodic points are dense, and  we can just choose any periodic point we like. In this case, to verify Assumptions \ref{A3} (v)-(vi) we just require regularity properties around the finite orbit of the periodic point.
\end{example}
\begin{example}
\label{ex:periodic-markov}
    Assumptions \ref{A3}  (i)-(iv)  hold when $(\Omega,\mathcal F,\mathbb P,\theta)$ is a Markov shift, such that the underlying Markov chain $\xi$ satisfies that $\mathbb P(\xi_i\in U_i; i\leq d)>0$ for all open sets $U_1,...,U_d$ in a neighborhood of a point $x_0$.
     In this case we can just take the fixed point  $\omega_0=(...,x_0,x_0,x_0,...)$ and require regularity properties around this point. For example, we can consider the case when $f_\omega$ and $R_\omega$ depend only on the $0$-th coordinates and are continuous with respect to it.
\end{example}

Under Assumption \ref{A3}, let $n_0$ be the period of $\omega_0$. Define  the operators $Q_t, t\in\mathbb R$ which map a function on $\cE_{\omega_0}$ to a function on $\cE_{\omega_0}$ by
$$
Q_tg=R_{\omega_0,it}\circ \cdots \circ R_{\theta^{n_0-1}\omega_0,it}g,
$$
where $R_{\omega,it}g=R_{\omega}(e^{itf_{\theta\omega}}g)$. 
\begin{theorem}\label{LLT rand}
Make Assumptions \ref{A2} and \ref{A3}. Suppose also that the number $\sigma$ from Remark \ref{rem} satisfies $\sigma>0$ (namely $\lim_{n\to\infty}\text{Var}(S_n^\omega)=\infty$, a.s). Assume also that the spectral radii of the operators  $Q_t, t\not=0$ are smaller than $1$. Then for $\mathbb P$-a.a. $\omega$ the sequence $S_n^\omega$ obeys the non-lattice local central limit theorem and the first-order Edgeworth expansions.   
\end{theorem}
The condition about the spectral radius of $Q_t$ is a standard assumption on the deterministic operators $Q_t$, see \cite{HH}. Since the operators $Q_t$ are quasi-compact, this condition means that the function $f(x_0,...,x_{n_0-1})=\sum_{j=0}^{n_0-1}f_{\theta^j\omega_0}(x_j)$ is aperiodic in the sense of \cite{HH}.

\begin{remark}
\label{rem:periodic-spectral-criterion} 
    A sufficient condition for the periodic spectral assumption in Theorem~\ref{LLT rand} can be obtained from the kernels at two successive points of the periodic orbit.
    Let $\omega_0$ be the periodic point of period $n_0$ appearing in that theorem, and write $\omega_j=\theta^j\omega_0$, with indices taken modulo $n_0$.
    Suppose that, for some index $s\in\{0,\ldots,n_0-1\}$, there are constants $\gamma_1,\gamma_2\in(0,1]$ and probability measures $\mu$ on $\cE_{\omega_{s+1}}$ and $\nu$ on $\cE_{\omega_{s+2}}$ such that $R_{\omega_s}(x,\cdot) \ge\gamma_1\mu(\cdot)$ for every $x\in\cE_{\omega_s}$ and $ R_{\omega_{s+1}}(y,\cdot) \ge\gamma_2\nu(\cdot)$ for every $y\in\cE_{\omega_{s+1}}$.
    Assume further that $(f_{\omega_{s+1}})_*\mu$ is non-lattice.
    Then the argument in Lemma~\ref{Key} below gives that
    \[
        \|R_{\omega_s,it}\circ R_{\omega_{s+1},it}\|_\infty
        \leq1-\gamma_1\gamma_2
          \left(1-\left|\int_{\cE_{\omega_{s+1}}}
                 e^{itf_{\omega_{s+1}}(y)}\,\mu(dy)\right|\right)
        <1,
        \qquad t\neq0.
    \]
    Since $R_{\omega_j}$ is a Markov kernel and $|e^{itf_{\omega_{j+1}}}|=1$, we have that each operator above has norm at most $1$.
    Since $Q^2_t = Q_t\circ Q_t$, the operators  $R_{\omega_s,it}$ and $R_{\omega_{s+1},it}$ occur consecutively and hence by submultiplicativity of norms we have that
     \[
        (r(Q_t))^2\leq\|Q_t^2\|_\infty<1
        \qquad\text{ for }\quad t\neq0.
    \]
\end{remark}

Next, we discuss Edgeworth expansions, where we focus on the case when $(\Omega,\cF,\mathbb P,\theta)$ is a subshift.
\begin{theorem}\label{Edge2}
Make Assumption \ref{A2} and suppose $\sigma>0$. Let $r\geq 2$.
Suppose 
\vskip0.1cm
\begin{enumerate}
    \item[(i)] $(\Omega,\cF,\mathbb P,\theta)$ is a topologically mixing subshift of finite type and that $\mathbb P$ is a Gibbs measure.
    \vskip0.1cm
    \item[(ii)]$\omega\to R_\omega$ and $\omega\to f_\omega$ depend only on the $0$-th coordinate.
    \vskip0.1cm
    \item[(iii)] (Diophanticity-see \cite{FL}) there is a periodic point $\omega_0\in\Omega$ and  $c,\delta>0$ and $m_0\in\mathbb N$ such that for all real $t$ with $|t|\geq\delta$,
$$
\|Q_t^{m_0}\|_\infty\leq 1-c|t|^{-\ell}
$$
for some $\ell$ satisfying $(r-1)\ell<2$.
\end{enumerate}
 Then for $\mathbb P$-a.a. $\omega$ the sequence $S_n^\omega$ obeys the Edgeworth expansions of order $r$.     
\end{theorem}

\section{Proof of Theorems \ref{LLT0}, \ref{LLT1} and \ref{LLT Red}}
Let $B_j$ denote the space of bounded measurable functions on $\cX_j$ with the supremum norm $\|\cdot\|_\infty$. For $z\in\mathbb C$ define $R_{j,z}:B_{j+1}\to B_j$ by
$$
R_{j,z}g(x)=\int e^{zf_{j+1}(y)}g(y)R_j(x,dy).
$$

We first need the following relatively simple result.
\begin{lemma}\label{Key}
Under Assumption \ref{A1},
for all  $t\in\mathbb R$  and $j\geq1$ we have 
$$
\|R_{j,it}\circ R_{j+1,it}\|_\infty\leq 1-\gamma^2\left(1-\left|\int e^{itf_{j+1}(y)}dm_{j+1}(y)\right|\right).
$$
\end{lemma}
\begin{proof}
Let $g:\cX_{j+2}\to\mathbb C$ be a measurable function with $\|g\|_\infty\leq1$. Take some $x\in\cX_{j}$.
Let us write $R_j(x,\cdot)=\gamma m_{j+1}(\cdot)+(1-\gamma)Q_j(x,\cdot)$ where $Q_j(x,\cdot)$ is a probability measure on $\cX_{j+1}$. Then, 
$$
R_{j}(x,dy_1)R_{j+1}(y_1,dy_2)=\gamma^2m_{j+1}(dy_1)m_{j+2}(dy_2)+(1-\gamma^2)P_j(x,dy_1,dy_2)
$$
where $P_j(x,\cdot,\cdot)$ is a probability measure on $\cX_{j+1}\times\cX_{j+2}$.
Now,   write 
$$
\left(R_{j,it}\circ R_{j+1,it}\right)g(x)=\iint e^{it(f_{j+1}(y_1)+f_{j+2}(y_2))}g(y_2)R_{j}(x,dy_1)R_{j+1}(y_1,dy_2)
$$
$$
=\gamma^2\left(\int e^{itf_{j+1}(y_1)}dm_{j+1}(y_1)\right) \left(\int e^{itf_{j+2}(y_2)}g(y_2)dm_{j+2}(y_2)\right)
$$
$$
+(1-\gamma^2)\int  e^{it(f_{j+1}(y_1)+f_{j+2}(y_2))}g(y_2)P_j(x,dy_1,dy_2).
$$
Thus, 
$$
\left|\left(R_{j,it}\circ R_{j+1,it}\right)g(x)\right|\leq \gamma^2\left|\int e^{itf_{j+1}(y_1)}dm_{j+1}(y_1)\right|+1-\gamma^2
$$
$$
=1-\gamma^2\left(1-\left|\int e^{itf_{j+1}(y_1)}dm_{j+1}(y_1)\right|\right).
$$
\end{proof}
\begin{remark}
Lemma \ref{Key} holds without the assumption that $f_j$ are bounded like in \cite[Lemma 11]{PeligradLLT2}, which dealt with a similar estimate but in the $L^1$ norm in the case when $m_j$ is the law of $X_j$.    
\end{remark}
\subsection{Proof of Theorem \ref{LLT0}}
Denote $\sigma_n=\sqrt{\text{Var}(S_n)}$.
Then by \cite[Proposition 3.2]{HW2}  and \cite[Corollary 28]{DHPTRF1} under $\phi$-mixing there are $\delta_0,c,C>0$ such that for all $t\in[-\delta_0,\delta_0]$ we have 
\begin{equation}\label{b}
\left|\mathbb E[e^{it S_n}]\right|\leq Ce^{-c\sigma_n^2 t^2}.    
\end{equation}
\begin{remark}[c.f. Remarks \ref{r1} and \ref{r2}]\label{bb}
  Without assuming uniform boundedness of $f_j$ we can still derive \eqref{b} by using \eqref{ppp} below, under Assumption \ref{A1} and the assumptions that $\sum_{j=1}^n\text{Var}_{m_j}(f_j)\geq c_1\sigma_n^2$ for $c_1>0$  and $\sup_j\|f_j\|_{L^{3}(m_j)}<\infty$. 
\end{remark}

Next, relying on \eqref{b} and arguing like in \cite[Chapter 2]{HK}, to prove the local limit theorem it is enough to show that for every compact $K\subset\mathbb R\setminus\{0\}$ we have 
\begin{equation}\label{Need}
\lim_{n\to\infty}\sigma_n\int_{K}\left|\mathbb E[e^{it S_n}]\right|dt=0.    
\end{equation}
Note that by \cite[Proposition 3.2]{HW2}  and \cite[Proposition 25]{DHPTRF1}, condition \eqref{up} is also sufficient for the first-order Edgeworth expansions.
Denote 
$$
D_n(t)=\sum_{j=1}^n d_j(t)
$$
where 
$$
d_j(t)=1-\left|\int e^{itf_j(y)}dm_j(y)\right|.
$$
Then by Lemma \ref{Key} and since $1-x\leq e^{-x}$ for all $x\geq0$ we have 
\begin{equation}\label{ppp}
|\mathbb E[e^{itS_n}]|\leq e^{\frac12\gamma^2}e^{-\frac12\gamma^2 D_n(t)}.    
\end{equation}
Indeed, by Assumption~\ref{A1}  for $n\ge 2$, we have that 
\[
    \|R_{n-1, it} \mathbf{1}\|_\infty \le 1 - \gamma d_n(t) \le e^{-\gamma^2 d_n(t)}.
\]
Then, apply Lemma~\ref{Key} to consecutive pairs with odd or even starting indices, using the estimate above to bound $R_{n-1,it}\mathbf{1}$ when the last operator $R_{n-1,it}$ is not included in a pair. 
Since one of the even-indexed or odd-indexed sums in $\sum_{j=2}^n d_j(t)$ is at least half the total sum and all remaining operators have norm at most one, we obtain \eqref{ppp}, with the use that $d_1(t) \le 1$.
Thus \eqref{Need} follows from the assumption of Theorem~\ref{LLT0}.

\subsection{Proof of Theorem \ref{LLT1}}
\subsubsection{General auxiliary results}
For every non-constant  random variable $Y$ and all real numbers $t$, we define
$$
d(Y,t)=\left\|e^{it(Y-Y')}-1\right\|_{L^2}=2\left\|\sin\left(\frac12t(Y-Y')\right)\right\|_{L^2}
$$
where $Y'$ is an independent copy of $Y$.
Next, set $\varphi_Y(t)=\bbE[e^{it Y}]$.

\begin{lemma}\label{LLL}
For every real $t$ we have 
$$
\frac14d^2(Y,t)\leq 1-|\varphi_Y(t)|\leq \frac{1}{2} d^2(Y,t).
$$
\end{lemma}
\begin{proof}
We have
$$
|\varphi_Y(t)|^2=1-2\bbE\left[\sin^2\left(\frac12t(Y-Y')\right)\right]=1-\frac12 d^2(Y,t)
$$ 
and so
$$
\frac14 d^2(X,t) =\frac{1-|\varphi_Y(t)|^2}{2}\leq 1-|\varphi_Y(t)|=\frac{1-|\varphi_Y(t)|^2}{1+|\varphi_Y(t)|}\leq 1-|\varphi_Y(t)|^2=\frac12 d^2(Y,t) .
$$
\end{proof}

The following lemma is close in spirit to \cite[Lemma 5.8]{DS}.
\begin{lemma}\label{Le3}
If $Y\in L^2$ then
for all $t,\eta\in\bbR$ we have 
$$
d^2(Y,t+\eta)\geq d^2(Y,t)+d^2(Y,\eta)-C|\eta|d(Y,t)\sqrt{\text{Var}(Y)},\,\,C=4\sqrt 2.
$$
If also $Y\in L^3$ then 
$$
d^2(Y,\eta)\geq 2\eta^2\text{Var}(Y)-\frac13|\eta|^3\bbE[|Y-Y'|^3]
$$
and so 
$$
d^2(Y,t+\eta)\geq  d^2(Y,t)+2\eta^2\text{Var}(Y)-\left(\frac13|\eta|^3\bbE[|Y-Y'|^3]+C|\eta|d(Y,t)\sqrt{\text{Var}(Y)}\right).
$$
\end{lemma}

\begin{proof}[Proof of Lemma \ref{Le3}]

Take $Z = Y - Y'$ and using $\cos(\alpha+\beta)=\cos\alpha\cos\beta-\sin\alpha\sin\beta$ we obtain that 
\begin{align*}
    d^2(Y,t+\eta)
        &= 2\mathbb E[1-\cos((t+\eta)Z)]\\
        &= 2\mathbb E[1-\cos(tZ)]  +2\mathbb E[1-\cos(\eta Z)] - 2\mathbb E[(1-\cos(tZ))(1-\cos(\eta Z))]\\
        &\qquad +2\mathbb E[\sin(tZ)\sin(\eta Z)]\\
        &= d^2(Y,t) + d^2(Y,\eta) + I_1 + I_2
\end{align*}
where
\[
    I_1=-2\mathbb E[(1-\cos(tZ))(1-\cos(\eta Z))],
    \qquad\text{ and }\qquad
    I_2=2\mathbb E[\sin(tZ)\sin(\eta Z)].
\]

Now to estimate $I_1$, we use that 
$$
2(1-\cos(tZ))=|e^{itZ}-1|^2
$$
and so by the Cauchy-Schwarz inequality,
$$
|I_1|^2\leq \bbE[|e^{itZ}-1|^4]\bbE[(1-\cos(\eta Z))^2]. 
$$
Since $|e^{itZ}-1|^2\leq4$, we have
\[
    \mathbb E[|e^{itZ}-1|^4] \leq4\mathbb E[|e^{itZ}-1|^2] =4d^2(Y,t).
\]
On the other hand, since $|1-\cos(\eta Z)|\leq |\eta Z|$ we see that 
\[
    \mathbb E[(1-\cos(\eta Z))^2] \leq\eta^2\mathbb E[Z^2] =2\eta^2\operatorname{Var}(Y).
\]
Consequently,
\[
    |I_1|^2 \le 8\eta^2d^2(Y,t) \operatorname{Var}(Y) \implies 
    |I_1| \le \sqrt{2} \left|\eta\right| d(Y,t) \sqrt{\operatorname{Var}(Y)}.
\]

Now for $I_2$ we have that 
$$
|\sin(tZ)|=\sqrt{1-\cos^2(tZ)}\leq 2\sqrt{2-2\cos(tZ)}=|e^{itZ}-1|.
$$
Thus, by the Cauchy-Schwarz inequality,
$$
|I_2|
= 2\left|\mathbb E\left[\sin(tZ)\sin(\eta Z)\right]\right|
\leq 2 \|\sin(tZ)\|_{L^2}\|\sin(\eta Z)\|_{L^2}
\le 2d(Y,t)\|\sin(\eta Z)\|_{L^2}.
$$
Using that $|\sin x|\leq |x|$ for all $x$ we see that 
$$
|I_2| \le 2 \sqrt{2} |\eta| d(Y,t) \sqrt{\operatorname{Var}(Y)}
$$

Combining these estimates, we see that 
\[
    d^2(Y, t+\eta) \ge d^2(Y,t) + d^2(Y,\eta) - 4\sqrt 2 |\eta|d(Y,t)\sqrt{\operatorname{Var}(Y)}
\]
Therefore, the first assertion in the lemma holds with $C  =4\sqrt{2}$.

As for the second assertion, suppose now that $Y \in L^3$, then Taylor's theorem yields that
\[
    \left|
        \cos u - \left(1 - \frac{u^2}{2}\right)
    \right|
    \le \frac{|u|^3}{6}
\]
for all $u\in\mathbb R$. 
Thus, in particular, 
\[
    2(1-\cos(u))\ge u^2 - \frac{|u|^3}{3}.
\]
We apply this with $u = \eta Z$ and taking expectations, we get that 
\[
    d^2(Y,\eta) 
        = 2 \mathbb E [1- \cos(\eta Z)]
            \ge \eta^2 \mathbb E[Z^2] 
            - \frac{|\eta|^3}{3} \mathbb E [|Z|^3]
        = 2\eta^2 \operatorname{Var}(Y) - \frac{|\eta|^3}{3} \mathbb E[|Y-Y'|^3].
\]

Finally, with $C = 4\sqrt{2}$, we obtain 
\[
    d^2(Y , t+\eta) \ge d^2(Y,t) + 2\eta^2\operatorname{Var}(Y) 
    -  
    \left(
        \frac{|\eta|^3}{3} \mathbb E[|Y-Y'|^3]
        + C|\eta|d(Y,t) 
        \sqrt{\operatorname{Var}(Y)}
    \right). 
\]
\end{proof}

\subsection{Proof of Theorem \ref{LLT1}}
As in the previous section, all that we need to do is to prove \eqref{Need}. 
First, take some compact set $K\subset\mathbb R\setminus\{0\}$. The following arguments are in the same spirit as the arguments in \cite[Section 5.2.4]{DS}.
Let us take some $\delta>0$ and let us cover $K$ by finitely many intervals $I_1,...,I_d$ whose length does not exceed $\delta$. By taking $\delta$ small enough we can ensure that $I_i\subset \mathbb R\setminus\{0\}$.
Then, it is enough to show that if $\delta$ is small enough then for all $i$,
$$
\lim_{n\to\infty}\sigma_n\int_{I_i}\left|\mathbb E[e^{itS_n}]\right|dt=0.
$$
Let us take some $i$ and denote $I_i=I$. Using \eqref{ppp}, it is enough to show that 
$$
\lim_{n\to\infty}\sigma_n\int_{I}e^{-\frac12\gamma^2D_n(t)}dt=0.
$$
Next, since $D_n(t)$ is increasing in $n$ and diverges pointwise to $\infty$ on $I$, we we see that it diverges to $\infty$ uniformly on $I$, that is
\begin{equation}\label{inf}
\lim_{n\to\infty}\inf_{t\in I}D_n(t)=\infty.    
\end{equation}
For $n$ large enough, let us take a point $a_n\in I$ such that $D_n(a_n)\leq \frac32\inf_{t\in I}D_n(t)$.
Now, take $Y_j=f_j$ on the probability space $(\cX_j,\cF_j, m_j)$. Let us view $(Y_j)_{j\geq0}$ as independent random variables on the product space $\prod(\cX_j,\cF_j, m_j)$. 
Then $\sup_j\|Y_j\|_{L^\infty}<\infty$. Next, by the Cauchy-Schwarz inequality, we have 
$$
\sum_{j=1}^n d(Y_j,t)\sqrt{\text{Var}(Y_j)}\leq \sqrt{\sum_{j=1}^n d^2(Y_j,t)}\sqrt{\sum_{j=1}^n\text{Var}(Y_j)}.
$$
Using this,  and then applying Lemma \ref{LLL} and
 Lemma \ref{Le3} we see that there are $\delta_0,c_0>0$ such that for all $n$ and all $\eta\in[-\delta_0,\delta_0]$,
$$
D_n(a_n+\eta)\geq \frac14\eta^2\sum_{j=1}^n\text{Var}(Y_j)-c_0 |\eta|\left(D_n(a_n)\sum_{j=1}^n\text{Var}(Y_j)\right)^{1/2}.
$$
Therefore, with $v_n=\sum_{j=1}^n\text{Var}(Y_j)$,
\begin{equation}\label{upppp}
    |\mathbb E[e^{i(a_n+\eta)S_n}]|\leq
    C e^{-\frac{\gamma^2}{8} \eta^2 v_n
        + c_0 |\eta|\sqrt{v_n D_n(a_n)}
    }.
\end{equation}
Recall our assumption that $v_n\geq c\text{Var}(S_n)=c\sigma_n^2$ for some $c>0$. Let us take some $R>0$. 
Then by \eqref{ppp}, we have 
\[
    \sigma_n \int_{t\in I,|t-a_n|\leq R\sqrt{\frac{D_n(a_n)}{v_n}}}
        |\mathbb E[e^{itS_n}]| \,dt
    \le CR (\sigma_n/\sqrt{v_n})\sqrt{D_n(a_n)} 
    e^{-\frac12\gamma^2\inf_{t\in I}D_n(t)}
    \to 0,
\]
where we used that $\inf_{t\in I}D_n(t)\geq 2D_n(a_n)/3$.
On the other hand, if $R$ is large enough, then there is $c_2>0$ such that
$$
\sigma_n\int_{t\in I, \,|t-a_n|\geq R\sqrt{\frac{D_n(a_n)}{v_n}}}|\mathbb E[e^{itS_n}]|dt\leq C_2\sigma_n\int_{-\infty}^\infty e^{-c_2\sqrt{v_nD_n(a_n)}|\eta|}d\eta=O\left(\frac{\sigma_n}{\sqrt{v_nD_n(a_n)}}\right)=o(1)
$$
where we used that $v_n\geq c\sigma_n^2$ and that $D_n(a_n)\to \infty$.

\subsection{Proof of Theorem \ref{LLT Red}}
We first need the following general fact.  For a random variable  $Y$ and $t\not=0$ denote 
\begin{equation}
\label{eq:D_y_t}
    D(Y,t)=\min_{\theta\in\mathbb R}\left\|\text{dist}(Y,\theta+(2\pi/t)\mathbb Z)\right\|_{L^2}.
\end{equation}
Recall that $d(Y,t)=2\|\sin(t(Y-Y')/2)\|_{L^2}$,
where $Y'$ is an independent copy of $Y$. Note that 
$$
|\mathbb E[e^{itY}]|^2=1-\frac12(d(Y,t))^2.
$$
\begin{lemma}
There is an absolute constant $C>0$ such that for every random variable $Y$ and all real $t\not=0$,
\begin{equation}\label{Fact0}
|t|D(Y,t)\leq Cd(Y,t).
\end{equation}
\end{lemma}
\begin{proof}
First, given a real number $y$, we have
$$
r=\text{dist}(y,\pi \bbZ)\leq \frac12\pi
$$
and so 
$$
|\sin y|=|\sin r|\geq C_1r=C_1\text{dist}(y,\pi \mathbb Z)
$$
where $C_1=\inf_{0< x\leq \frac12\pi}\frac{|\sin x|}{x}$.

Now, let $Y'$ be an independent copy of $Y$. Then for every realization of $y=Y'$ of $Y'$ we have 
$$
(D(Y,t))^2\leq \left\|\text{dist}\left(Y, y'+\left(\frac{2\pi}{t}\right)\bbZ\right)\right\|_{L^2}^2.
$$
Thus, by integrating with respect to the law of $Y'$, we see that 
$$
D(Y,t)\leq \left\|\text{dist}\left(Y-Y', \left(\frac{2\pi}{t}\right)\bbZ\right)\right\|_{L^2}:=\cD(Y-Y',t).
$$
Now, notice that 
$$
\cD(Y-Y',t)=\frac{2}{|t|}\left\|\text{dist}\left(\frac{t(Y-Y')}{2},\pi\bbZ\right)\right\|_{L^2}.
$$
By applying the previous estimate with $y=t(Y-Y')/2$ we conclude that 
$$
\cD(Y-Y',t)\leq \frac{2C_1^{-1}}{|t|}\left\|\sin\left(\frac12t(Y-Y')\right)\right\|_{L^2}=  \frac{1}{|t|}C_1^{-1}d(Y,t).
$$
\end{proof}

Now we are ready to complete the proof of Theorem \ref{LLT Red}.
Notice first that
$$1-\left|\int e^{itf_j(y)}dm_j(y)\right|^2=2\int\int\sin^2\left(\frac{t}{2}(f_j(x)-f_j(y))\right)dm_j(x)dm_j(y).
$$
Suppose that 
 for some $t\not=0$,
$$
\sum_{j=1}^{\infty}\left(1-\left|\int e^{itf_j(y)}dm_j(y)\right|\right)<\infty.
$$
Then since $2\left(1-\left|\int e^{itf_j(y)}dm_j(y)\right|\right)\geq 1-\left|\int e^{itf_j(y)}dm_j(y)\right|^2$
we have
$$
\sum_{j=1}^{\infty}\left(1-\left|\int e^{itf_j(y)}dm_j(y)\right|^2\right)<\infty.
$$
Take $Y_j(x)=f_j(x)$ on the probability space $(\cX_j,\mathcal F_j,m_j)$, and assume that the $Y_j$ are independent. Take $\theta_j\in\mathbb R$ that minimizes in the definition of $D(Y_j,t)$. Let $Z_j$ be an integer-valued random variable satisfying 
$$
\text{dist}(Y_j,\theta_j+(2\pi/t)\mathbb Z)=|Y_j-\theta_j-(2\pi/t)Z_j|.
$$
Let 
$$
r_j=Y_j-\theta_j-(2\pi/t)Z_j.
$$
Then 
$$
\int r_j^2\,dm_j\leq D^2(f_j,t)
$$
and so 
$$
\sum_j\int r_j^2\,dm_j<\infty.
$$
This gives us \eqref{Red0} with $a_j=\theta_j+\int r_j\,dm_j$.
Relying on \eqref{Fact0}, and using Kolmogorov's three series theorem with the random variables $r_j-\int r_j\,dm_j$ 
 we obtain \eqref{Red}.  
\qed

\subsection{Proof of Theorem~\ref{Edge1}}
Theorem \ref{Edge1} follows from \eqref{ppp} together with \cite[Proposition 25]{DHPTRF1}.

\subsection{Proof of Proposition~\ref{EdgeProp}}
Condition \eqref{EdgeDec} clearly holds when $A\geq \delta_0$. To handle the case of smaller $A$'s, it is enough to show that 
for all compact $K\subset\mathbb R\setminus\{0\}$ we have 
\begin{equation}\label{show}
\limsup_{j\to\infty}\sup_{t\in K}\left|\int e^{itf_j(y)}dm_j(y)\right|<1.    
\end{equation}
Indeed, under \eqref{show} we can use \eqref{ppp} to conclude that for all $n\in\mathbb N$,
$$
\int_{K}|\mathbb E[e^{itS_n}]|dt\leq Ce^{-cn}
$$
for some constants $C,c>0$ that may depend on $K$ but not on $n$.

Next, let us prove \eqref{show}. If \eqref{show} fails, then, possibly along a subsequence, there are $t_j\in K$ such that 
$$
\lim_{j\to\infty}\left|\int e^{itf_j(y)}dm_j(y)\right|=1.
$$
Now, since
$$
\left|\int e^{itf_j(y)}dm_j(y)\right|^2=1-2\iint \sin^2(t_j(f_j(x)-f_j(y))/2)dm_j(x)dm_j(y)
$$
we get that
\begin{equation}\label{uu}
\lim_{j\to\infty}\iint \sin^2(t_j(f_j(x)-f_j(y))/2)dm_j(x)dm_j(y)=0.  
\end{equation}
Next, for all integer $s\geq2$ and all $\alpha\in\mathbb R$ we have
$$
\sin(s\alpha)=\sum_{k=1}^s\binom{s}{k}\cos^{s-k}(\alpha)\sin^k(\alpha)\sin(\pi k/2).
$$
Therefore, there is $c_s>0$  depending only on $s$ such that for all $\alpha\in\mathbb R$, 
$$
\sin^2(s\alpha)\leq c_s\sin^2(\alpha). 
$$
This, together with \eqref{uu} implies that for all integers $s$ we have 
$$
\lim_{j\to\infty}\iint \sin^2(st_j(f_j(x)-f_j(y))/2)dm_j(x)dm_j(y)=0.
$$
Using that 
$$
\left|\int e^{ist_jf_j(y)}dm_j(y)\right|^2=1-2\iint \sin^2(s t_j(f_j(x)-f_j(y))/2)dm_j(x)dm_j(y)
$$
we conclude that 
\begin{equation}\label{c}
 \lim_{j\to\infty}\left|\int e^{i st_jf_j(y)}dm_j(y)\right|=1.   
\end{equation}
However, take $s$  large enough such that $\inf_{t\in K}|st|>\delta_0$. Then   $\inf_j|st_j|>\delta_0$ and so
$$
\left|\int e^{i st_jf_j(y)}dm_j(y)\right|\leq 1-c|st_j|^{-\ell}
$$
which together with the fact that $\sup_j|st_j|\leq\sup_{t\in k}|st|<\infty$ contradicts \eqref{c}.
\qed

\section{Proof of Theorems \ref{LLT rand0} and \ref{LLT rand}}
\subsection{Proof of Theorem \ref{LLT rand0}}
\,
\vskip0.1cm
\paragraph{\textbf{Proof of Theorem \ref{LLT rand0} (i)}} First, since $v(\omega):=\text{Var}_{m_\omega}(f_\omega)$ is not zero $\mathbb P$-a.s. by the mean ergodic theorem we see that 
$$
\lim_{n\to\infty}\frac1n\sum_{j=1}^nv(\theta^j\omega)=\int v(\omega')d\mathbb P(\omega')>0
$$
for $\mathbb P$-a.a $\omega$. Thus, $\mathbb P$-a.s. the condition in Theorem \ref{LLT1} that for all $n$ large enough we have 
$$
\sum_{j=1}^nv(\theta^j\omega)\geq c_\omega\sum_{j=1}^n\text{Var}(f_{\theta^j\omega}(X_{\omega,j}))
$$
where $c_\omega>0$,  is in force. 

To verify the second condition in Theorem \ref{LLT1} we need to show that $\mathbb P$-a.s. for all nonzero $t\in\mathbb R$ we have 
\begin{equation}\label{div1}
\sum_{j=1}^{\infty}\left(1-\left|\int e^{itf_{\theta^j\omega}(y)}dm_{\theta^j\omega}(y)\right|\right)=\infty.    
\end{equation}
Fix some $t\in\mathbb R\setminus\{0\}$. Denote 
$$
\nu_\omega(t)=1-\left|\int e^{itf_{\omega}(y)}dm_{\omega}(y)\right|.
$$
Let $h=\frac{2\pi}{t}$. Let $\Delta_h=\Omega\setminus\Gamma_h$, where $\Gamma_h=\{\omega\in\Omega: |m_\omega(e^{it f_\omega})|=1\}$.
Then by assumption $\mathbb P(\Delta_h)>0$ and so for all nonzero $t$ we have 
$$
\nu(t):=\int \nu_\omega(t)d\mathbb P(\omega)>0.
$$
Now, by the mean ergodic theorem for all $t\in\mathbb R\setminus\{0\}$ for $\mathbb P$-a.a. $\omega$ we have
\begin{equation}\label{div2} 
\lim_{n\to\infty}\frac1n\sum_{j=1}^n\nu_{\theta^j\omega}(t)=\nu(t).    
\end{equation}
This shows that \eqref{div1} holds $\mathbb P$-a.s. for a given $t$. Thus, \eqref{div1} holds true $\mathbb P$-a.s. for all $t\in\mathbb Q\setminus\{0\}$. Now, let us take some $\omega$ such that \eqref{div2} holds for  all $t\in\mathbb Q\setminus\{0\}$. 
 Let $K=\text{ess-sup}\|f_\omega\|_{\infty}<\infty$. 
 Given an irrational $s\in\mathbb R\setminus\{0\}$, and $t\in\mathbb Q\setminus\{0\}$ for all $\omega\in\Omega$ we have
$$
\sup_n\left|\frac1n\sum_{j=1}^n\nu_{\theta^j\omega}(s)-\frac1n\sum_{j=1}^n\nu_{\theta^j\omega}(t)\right|\leq K|t-s|
$$
and also 
$$
|\nu(s)-\nu(t)|\leq K|t-s|.
$$
Thus, if also $K|t-s|<\frac14\nu(s)$ then for all $n$ and $\omega\in\Omega$,
$$
\frac1n\sum_{j=1}^n\nu_{\theta^j\omega}(s)\geq \frac1n\sum_{j=1}^n\nu_{\theta^j\omega}(t)-K|t-s|.
$$
Taking the limit as $n\to\infty$ we see that 
$$
\liminf_{n\to\infty}\frac1n\sum_{j=1}^n\nu_{\theta^j\omega}(s)\geq \nu(t)-K|t-s|\geq \nu(s)-2K|t-s|>\frac12\nu(s)>0.
$$
We thus conclude that $\mathbb P$-a.s. \eqref{div1} holds for every nonzero $t$. This completes the proof of Theorem \ref{LLT rand0} (i).
\qed
\vskip0.1cm

\paragraph{\textbf{Proof of  Theorem \ref{LLT rand0} (ii)}} 
For this purpose, we will verify the assumption of Theorem \ref{Edge1} with the sequence $X_j=X_{\omega,j}$ for almost all $\omega$. 

 First,  denote by $B_0$ the set of all $\omega\in\Omega$ such that $|m_{\omega}(e^{itf_\omega})|\leq 1-c_\omega|t|^{-\ell}$ for all $t$ such that $|t|\geq \delta_\omega$. Then for all $\omega\in B_0$  for all $t\in\mathbb R\setminus\{0\}$ we have 
 $$
|m_{\omega}(e^{itf_\omega})|<1.
 $$
 Indeed, if $|m_{\omega}(e^{itf_\omega})|=1$ for some $t\not=0$ then $f_\omega$ takes values on a set of the form $a_\omega+(\frac{2\pi}t)\mathbb Z$ and so $|m_{\omega}(e^{isf_\omega})|=1$ for all multiples of $t$, which contradicts the definition of $B_0$. 
Let $K\subset\mathbb R\setminus\{0\}$ be a compact set. Then for $\omega\in B_0$ we have 
$$
c_K(\omega):=\sup_{t\in K}\left|\int e^{itf_\omega (y)}dm_\omega(y)\right|<1.
$$
For $0<b<1$ let $B=B_K=\{\omega\in B_0\in\Gamma: c_K(\omega)<b\}$. Then by taking $b=b_K$ close enough to $1$ we get that $\mathbb P(B)>0$. Consider the sequence of visiting times to $B$, that is $0\leq m_1(\omega)<m_2(\omega)<...$ such that for all $n\geq0$ we have $\theta^m\omega\in B$ if and only if $m=m_i(\omega)$ for some $i$. By the Poincar\'{e} recurrence theorem, this sequence is almost surely infinite.
Let $g(\omega)=\mathbb I(\omega\in B)$. Set $M_n(\omega)=\left|\{i: m_i(\omega)\leq n\}\right|$. Then 
$$
M_n(\omega)=\sum_{j=0}^{n}g(\theta^j\omega)
$$
and so by the mean ergodic theorem
$$
\lim_{n\to\infty}\frac{M_n(\omega)}{n}=\mathbb P(B),\text{ a.s.}
$$
Using this we conclude that $\mathbb P$-a.s.  for all $n$ large enough
\begin{equation}\label{cpt}
\sup_{t\in K}\exp\left(-\frac12\gamma^2\sum_{m=1}^n\left(1-\left|\int e^{itf_{\theta^m\omega}(y)}dm_{\theta^m\omega}(y)\right|\right)\right)    
\end{equation}
$$
\leq \sup_{t\in K}\exp\left(-\frac12\gamma^2M_n(\omega)(1-b)\right)
\leq \exp\left(-\frac14\gamma^2n\mathbb P(B)(1-b)\right).
$$
We thus see that the product decays exponentially fast to $0$.

Using \eqref{cpt}, what remains in order to prove Theorem \ref{LLT rand0} (ii) is to show that for all $A>0$ large enough and all $B>0$ for $\mathbb P$-a.a. $\omega$ we have 
\begin{equation}\label{long}
\lim_{n\to\infty}n^{r/2}\int_{A\leq |t|\leq Bn^{(r-1)/2}}\exp\left(-\frac12\gamma^2\sum_{m=1}^n\left(1-\left|\int e^{itf_{\theta^m\omega}(y)}dm_{\theta^m\omega}(y)\right|\right)\right)dt=0.    
\end{equation}
Let us take $\delta>0$ large enough and $c>0$ small enough such that the set  
$$
B=\left\{\omega\in B_0: \delta_\omega\leq \delta, c_\omega>c\right\}
$$
has positive probability. Let $m_1(\omega)<m_2(\omega)<...$ be the visiting times to the set $B$. Then for every real $t$ with $|t|\geq\delta$ we have 
$$
\exp\left(-\frac12\gamma^2\sum_{m=1}^n\left(1-\left|\int e^{itf_{\theta^m\omega}(y)}dm_{\theta^m\omega}(y)\right|\right)\right)\leq \exp\left(-\frac12c\gamma^2M_n(\omega)|t|^{-\ell}\right)
$$
where $M_n(\omega)$ is the number of visits until time $n$, which is a.s. at least of order $n\mathbb P(B)/2$. Thus there exists $c_0>0$ such that for all real $t$ with $|t|\geq\delta$ we have 
$$
\exp\left(-\frac12\gamma^2\sum_{m=1}^n\left(1-\left|\int e^{itf_{\theta^m\omega}(y)}dm_{\theta^m\omega}(y)\right|\right)\right)\leq \exp\left(-c_0 n^{1-\frac12(r-1)\ell}\right). 
$$
Therefore \eqref{long} holds, and the proof of Theorem \ref{LLT rand0} (ii) is complete.
\qed

\subsection{Proof of Theorem \ref{LLT rand}}
Denote 
$$
R_{\omega,n}=R_{\omega}\circ R_{\theta\omega}\circ\cdots\circ R_{\theta^{n-1}\omega}.
$$
Then by \cite[Proposition 3.2]{HW2} and \cite[Corollary 28]{DHPTRF1} applied with the sequence of operators $R_j=R_{\theta^j\omega}$ and the fact that $\text{Var}(S_n^\omega)$ grows linearly fast in $n$ (see \cite[Theorem 2.1]{Kifer1998}), we see that $S_n^\omega$ obeys the central limit theorem and there exists $c_0,C_0>0$ such that for $\mathbb P$-a.a. $\omega$ for all $|t|\leq c_0$ and all $n$ large enough we have
$$
|\mathbb E[e^{it S_n^\omega}]|\leq C_0e^{-cnt^2}.
$$
Thus, arguing like in \cite[Chapter 2]{HK}, to complete the proof of the local central limit theorem it is enough to show that for all compact $K\subset\mathbb R\setminus\{0\}$ we have 
\begin{equation}\label{up}
 \lim_{n\to\infty}\sqrt n\sup_{t\in K}|\mathbb E[e^{it S_n^\omega}]|=0.   
\end{equation}
Note that by \cite[Proposition 3.2]{HW2}  and \cite[Proposition 25]{DHPTRF1}\, \eqref{up} is also sufficient for the first-order Edgeworth expansions.
To prove \eqref{up}, we will use the properties in Assumption \ref{A2}. First, by upper semi-continuity of the spectral radius and the arguments in the proof of \cite[Corollary III.13]{HH}, there exist $\epsilon=\epsilon_K>0$ and $C=C_K>0$ such that for all $n\in\mathbb N$ we have
$$
\sup_{t\in K}\|Q_t^n\|_\infty\leq C\left(1-\epsilon\right)^n.
$$
Let us take  $m_0\in\mathbb N$ such that $C(1-\epsilon)^{m_0}<1/4$.
Let us take a neighborhood $U$ of $\omega_0$ 
such that
$$
\sup_{\omega\in U}\sup_{t\in K}\|R_{\omega,it}^{m_0n_0}-Q_t^{m_0}\|_\infty<1/4.
$$
Then, 
$$
\sup_{\omega\in U}\sup_{t\in K}\|R_{\omega, it}^{m_0n_0}\|_\infty<\frac12.
$$
Now, let us take the visiting times $k_1(\omega)<k_2(\omega)<...$ to $U$ and let $\ell_{s}(\omega)=k_{sm_0n_0}(\omega)$. Since $\|R_{\omega,t}\|_\infty\leq1$ we thus conclude that 
$$
\sup_{t\in K}\|R_{\omega,it}^n\|_\infty\leq 2^{-\left|\{i\leq n: \ell_i(\omega)\leq n\right|}.
$$
Since 
$$
\lim_{n\to\infty}\frac{\left|\{i\leq n: \ell_i(\omega)\leq n\}\right|}{n}=\frac{\mathbb P(U)}{m_0n_0}, \,\,\,\,\mathbb P-a.s.
$$
we see that there is $c_K>0$ such that for $\mathbb P$-a.a. $\omega$ for all $n$ large enough we have
$$
\sup_{t\in K}\|R_{\omega,it}^n\|_\infty\leq e^{-c_Kn}.
$$
Thus \eqref{up} holds with an exponential rate of decay.

\subsection{Proof of Theorem \ref{Edge2}}
Let $\omega_0=(a_j)_{j\in\mathbb Z}$ be the periodic point from the theorem and let $n_0$ be its period. Let us take $|t|\geq\delta$. Then by the mean ergodic theorem for $\mathbb P$-a.a. $\omega=(b_j)_{j\in\mathbb Z}$ there exist $c_0>0$ and $n_\omega$ such that for all $n\geq n_\omega$ we have 
$$
\sum_{j=0}^{\left[\frac{n}{n_0m_0}\right]}\mathbb I(b_j=a_j;\, \forall jn_0m_0\leq j< (j+1)n_0m_0)\geq c_0 n.
$$
Here $\mathbb I(\cdot)$ equals $1$ when the formula in the brackets holds, and $0$ otherwise. 
Notice that as $R_{\omega,t}$ depends only on the zero coordinate $b_0$, we have that on the above set and for the indices $j$ such that $b_j=a_j$ for all $jn_0m_0\leq j<(j+1)n_0m_0$, we have
$$
R_{\theta^{jn_0m_0+1}\omega, it}^{n_0m_0}=Q_t^{m_0}.
$$
Now, as  $\|R_{\omega,it}\|_\infty\leq1$, by submultiplicativity of operator norms for all $n\geq n_\omega$ and $|t|\geq\delta$ we have
$$
\left\|R_{\omega,it}^n\right\|_\infty\leq \left(1-c|t|^{-\ell}\right)^{c_0 n}\leq e^{-c_0n|t|^{-\ell}}.
$$
If also $|t|\le Bn^{(r-1)/2}$, this is at most $e^{-cc_0B^{-\ell}n^\epsilon}$, where $\epsilon=1-\frac12(r-1)\ell>0$.
Now the result follows from \cite[Proposition 3.2]{HW2}  and \cite[Proposition 25]{DHPTRF1}.
 \qed


\section{Application to Disordered Quantum Trajectories}

In this section, we explain how the preceding results apply to posterior-state observables of a disordered quantum trajectory.
We remind the reader of the mathematical setting discussed in \cite{traj_paper} briefly below. 

Quantum trajectories describe a microscopic evolution of the quantum states of a system that undergoes repeated measurements.
In a more realistic setting, the repeated measurements themselves may be noisy, resulting in a classically disordered quantum trajectory. 
We assume that the system which undergoes repeated measurements has a finite-dimensional Hilbert space $\mathbb C^d$.
Thus, the quantum states of the system are represented by density matrices 
\[
    \mathbb S_d=\{\rho\in\mathbb M_d(\mathbb C):
        \rho\geq0,\ \operatorname{Tr}\rho=1\}.
\]
Here $\mathbb M_d(\mathbb C)$ denotes the space of $d\times d$ complex matrices. 
We equip $\mathbb S_d$ with the Borel $\sigma$-algebra of the usual matrix topology. 

We model the classical disorder by an invertible ergodic probability-preserving system $(\Omega,\cF,\mathbb P,\theta)$.
We assume that the outcome space, $(\Xi,\cX)$, is a standard Borel space.
We start with a random instrument: $\omega\mapsto\mathfrak I_\omega$: for each $\omega$, the map $A\mapsto\mathfrak I_\omega(A)$ is countably additive, where each $\mathfrak I_\omega(A)$ is a completely positive linear map on $\mathbb M_d(\mathbb C)$ with the condition that $\mathfrak I_\omega(\Xi)$ is trace-preserving. 
For each $A\in\cX$ and $X\in\mathbb M_d(\mathbb C)$, the map $\omega\mapsto\mathfrak I_\omega(A)(X)$ is measurable.
We call the ergodic probability-preserving system $(\Omega,\cF,\mathbb P,\theta)$ the disorder environment, and this dictates which instrument process is applied at each time by translations via $\theta$: for each disorder realization $\omega\in\Omega$ at time $j$ the measurement is done via the instrument $\mathfrak I_j(\omega):=\mathfrak I_{\theta^j\omega}$. 

By \cite[Proposition~2.2]{traj_paper}, after modifying the instrument family on a $\mathbb P$-null set if necessary, the random instrument admits a measurable density realization
\[
    \mathfrak I_\omega(A)(X)
        =\int_A\cT_{x;\omega}(X)\,\nu_\omega(dx),\qquad A\in\cX,
\]
where $\nu_\omega$ is a measurable probability kernel used as a reference measure.
The map $(\omega,x)\mapsto\cT_{x;\omega}$ is (jointly) measurable, each $\cT_{x;\omega}$ is completely positive, and the non-selective map $\Phi_\omega=\mathfrak I_\omega(\Xi)$ is trace preserving.

For an incoming state $\rho$, the matrix $\mathfrak I_\omega(A)(\rho)$ is the unnormalized state associated with observing an outcome in $A$; its trace is the probability of this event. In particular, the next outcome has distribution
\[
    \Gamma_\omega(\rho,dx)
      =\operatorname{Tr}(\cT_{x;\omega}(\rho))\,\nu_\omega(dx).
\]
Conditioned on the observed outcome $x$, the posterior state is
\[
    F_\omega(\rho,x)
        =\frac{\cT_{x;\omega}(\rho)}
            {\operatorname{Tr}(\cT_{x;\omega}(\rho))}
\]
whenever the denominator above is non-zero. 
On the zero-denominator set we define $F_\omega(\rho,x)$ to be a fixed state $\rho_*\in\mathbb S_d$ and this set has zero $\Gamma_\omega(\rho, \,\cdot\,)$-probability. 

More precisely, in this setting if $\vartheta:\Omega\to\mathbb S_d$ is a random initial state and if $\bar x=(x_1, x_2, \ldots )\in \Xi^{\mathbb N}$ is a sequence of successive measurement outcomes, the corresponding trajectory started at $\vartheta$ in the disorder realization $\omega$ is defined recursively by
\[
    \rho_0^{\omega,\vartheta}=\vartheta(\omega),
    \qquad\text{ and }\qquad
    \rho_{j+1}^{\omega,\vartheta}
        =F_{\theta^j\omega}(\rho_j^{\omega,\vartheta},x_{j+1}),
    \quad \text{ for } \  j\geq0.
\]
Thus $(\rho_j^{\omega,\vartheta})_{j\geq0}$ is a measurable function of the outcome sequence $\bar x$, with both the environment $\omega$ and the initial state $\vartheta(\omega)$ fixed.

For this initial state, let $\mathbb Q_{\vartheta;\omega}$ denote the quenched law of the measurement outcomes on $(\Xi^{\mathbb N}, \cX^{\otimes\mathbb N})$.
It is determined by the cylinder probabilities
\[
    \mathbb Q_{\vartheta;\omega}\{x_1\in A_1,\cdots , x_n \in A_n\}
    =
    \operatorname{Tr}\!\left[
        \bigl(\mathfrak I_{\theta^{n-1}\omega}(A_n)\circ\cdots\circ
        \mathfrak I_\omega(A_1)\bigr)(\vartheta(\omega))
    \right],
\]
for $ A_1,\ldots,A_n\in\cX$. 
We refer the reader to \cite{QT_Pur} for the detailed description of this disordered setting and construction of the quenched law and to \cite{QT_LLN} for ergodic properties of the quenched law. 

Under $\mathbb Q_{\vartheta;\omega}$, conditioned on $x_1,\ldots,x_j$, the next outcome $x_{j+1}$ has distribution $\Gamma_{\theta^j\omega}(\rho_j^{\omega,\vartheta},\,\cdot\,)$.
Hence the conditional distribution of the next posterior state depends on the past only through the current posterior state and the fixed disorder realization.
This means that the trajectory of an initial state is an inhomogeneous Markov chain on $\mathbb S_d$, with one-step posterior kernel given by 
\begin{equation}
\label{eq:qt-posterior-kernel}
    K_\omega(\rho,B)
      =\int_\Xi 1_B(F_\omega(\rho,x))
           \operatorname{Tr}(\cT_{x;\omega}(\rho))\,\nu_\omega(dx),
            \qquad B\in\cB(\mathbb S_d).
\end{equation}

The law of the entire posterior-state trajectory $(\rho_j^{\omega,\vartheta})_{j\geq0}$ on $\mathbb S_d^{\mathbb N_0}$ is denoted by $\mathbb K_{\vartheta;\omega}$.
Equivalently, it is the pushforward
\begin{equation}
\label{eq:qt-state-path-law}
 \mathbb K_{\vartheta;\omega}
    =
    \left(
        \bar x \mapsto (\rho_j^{\omega,\vartheta}(\bar x))_{j\geq0}
    \right)_* \mathbb Q_{\vartheta;\omega}.
\end{equation}
Its initial distribution is $\delta_{\vartheta(\omega)}$, and its successive transition kernels are  $K_{\theta^j\omega}$, for each $j\geq0$.
Both $\mathbb Q_{\vartheta;\omega}$ and $\mathbb K_{\vartheta;\omega}$ depend on the initial state $\vartheta(\omega)$.
We use the same notation $\rho_j^{\omega,\vartheta}$ for the posterior process on the outcome space and the coordinate process under $\mathbb K_{\vartheta;\omega}$.

\subsection{Local limits for posterior-state sums}

We first apply Theorem~\ref{LLT rand0} to the disordered trajectories described above.
Assume that the following one-step Doeblin-type condition is satisfied by the kernel in \eqref{eq:qt-posterior-kernel}.
Suppose that there is a measurable family $\eta_\omega$ of probability measures on $\mathbb S_d$ and a deterministic constant $\varepsilon\in(0,1)$ such that
\begin{equation}\label{eq:qr-one-step}
    K_\omega(\rho,\cdot)\geq\varepsilon\eta_\omega(\cdot)
\end{equation}
for $\mathbb P$-almost every $\omega$, simultaneously for all $\rho\in\mathbb S_d$.
Examples satisfying this one-step minorization are given in \cite[Appendix~A]{traj_paper}.
We fix such a minorizing family $\eta_\omega$ for the rest of this subsection.

Let $f:\Omega\times\mathbb S_d\to\mathbb R$ be a jointly measurable bounded real-valued function, and write $f_\omega(\rho):=f(\omega,\rho)$.
Physically relevant choices include, but are not limited to, 
(i) the expectation $\operatorname{Tr}(H_\omega\rho)$ of a quantum observable $H_\omega$ in state $\rho$;
(ii) the purity $\operatorname{Tr}(\rho^2)$; and
(iii) the von Neumann entropy $-\operatorname{Tr}(\rho\log\rho)$, with $0\log0:=0$.
Here $H_\omega$ is a measurable, uniformly bounded family of self-adjoint matrices. 
All these functions are continuous in $\rho$ and uniformly bounded in the fixed finite dimension.
Their non-lattice behavior must be checked under the chosen minorizing measure.

For a measurable random initial state $\vartheta:\Omega\to\mathbb S_d$, define
\begin{equation}
\label{eq:qt-sum}
    S_n^{\omega,\vartheta}
      :=\sum_{j=1}^n f_{\theta^j\omega}
                           (\rho_j^{\omega,\vartheta}),
    \qquad\text{ and}\qquad 
    V_n^{\omega,\vartheta}
      :=\operatorname{Var}_{\mathbb K_{\vartheta;\omega}}
                           (S_n^{\omega,\vartheta}).
\end{equation}
Under a uniform posterior Doeblin condition and positive asymptotic variance, a quenched central limit theorem and standardized Berry--Esseen bounds for these sums were obtained in \cite[Theorem~3]{traj_paper}.
We will strengthen these conclusions to a non-lattice local central limit theorem and a first-order Edgeworth expansion below.
To this end, we identify
\begin{equation}\label{eq:qt-dictionary}
    \cE_\omega=\mathbb S_d,\qquad
    X_{\omega,j}=\rho_j^{\omega,\vartheta},\qquad
    R_\omega=K_\omega,\qquad
    \gamma=\varepsilon,\qquad
    m_\omega=\eta_{\theta^{-1}\omega},
\end{equation}
where the quantities on the left are as in Theorem~\ref{LLT rand0} and those on the right are defined above.
Then \eqref{eq:qr-one-step} is precisely the minorization
$R_\omega(\rho,\cdot)\geq\gamma m_{\theta\omega}(\cdot)$ required in that theorem.
For almost every fixed environment, it also gives Assumption~\ref{A1} with $m_j=\eta_{\theta^{j-1}\omega}$.
The forward estimates in Theorem~\ref{LLT rand0} apply to the prescribed initial state.

It is convenient to introduce the scalar probability measures 
\begin{equation}
\label{eq:qt-scalar-law}
    \mu_\om = (f_{\theta\om})_* \eta_\omega . 
\end{equation}
Equivalently, $\mu_\om(B) = \eta_\omega \{\rho\in\mathbb S_d: f_{\theta\om}(\rho) \in B\}$ for all Borel sets $B\subset\mathbb R$. 
Now we denote 
\[
    \widehat\mu_\om(t) := \int_{\mathbb{R}} e^{itu} \mu_\om(d u).
\]
Since $(f_\om)_* m_\om = \mu(\theta^{-1}\om)$ and as $\theta$ is $\mathbb P$-preserving, the condition in Theorem~\ref{LLT rand0} is equivalent to
\begin{equation}
\label{eq:qt-nonlattice}
    \mathbb P \{\om : |\widehat\mu_\om(t)| < 1\} > 0, \qquad 
    \text{for each $t\in\mathbb R\setminus \{0\}$}. 
\end{equation}
Equivalently, for every fixed lattice spacing $h>0$, there is an event $B_h\in \mathcal F$ such that $\mathbb P(B_h) >0$ with 
\[
    \mu_\om(a+h\mathbb Z) < 1
\]
for each $\om\in B_h$ and for every $a\in\mathbb R$. 
Therefore, once a minorizing family of measures $\eta_\omega$ on $\mathbb S_d$ satisfying \eqref{eq:qr-one-step} has been chosen,  we may apply Theorem~\ref{LLT rand0} to jointly measurable and uniformly bounded functions $f:\Omega\times\mathbb S_d\to\mathbb R$ that satisfy the condition \eqref{eq:qt-nonlattice} to obtain the following corollary.

\begin{corollary}[Local limits for posterior-state sums]
\label{cor:llt_qt_1}
    Suppose \eqref{eq:qr-one-step} and \eqref{eq:qt-nonlattice} hold.
    Then there are deterministic constants $c_*, C_* >0 $ such that, for every measurable initial state $\vartheta:\Om\to\mathbb S_d$ and $\mathbb P$-almost every $\omega$ we have that
    \[
        c_*n \le V_n^{\omega,\vartheta} \le C_*n 
    \]
    for all sufficiently large $n$. 
    In particular $\operatorname{Var}_{\mathbb K_{\vartheta;\omega}}(S_n^{\omega,\vartheta}) \to \infty$ and the sums $S_n^{\omega,\vartheta}$ satisfy the non-lattice local central limit theorem under $\mathbb K_{\vartheta;\omega}$ as defined in Definition~\ref{dfn:LCLT}.
    Furthermore, they also satisfy the first-order Edgeworth expansion; explicitly,
    \begin{equation}
    \label{eq:qt-first-edgeworth}
    \begin{aligned}
    &\sup_{t\in\mathbb R}
    \left|
        \mathbb K_{\vartheta;\omega}
        \left\{
            \frac{S_n^{\omega,\vartheta}
                -\mathbb E_{\mathbb K_{\vartheta;\omega}}
                    [S_n^{\omega,\vartheta}]}
                 {\sqrt{V_n^{\omega,\vartheta}}}
            \leq t
        \right\}
        -\Phi(t)
        -\frac{\kappa_{3,n}^{\omega,\vartheta}}
                {6(V_n^{\omega,\vartheta})^{3/2}}
            (1-t^2)\varphi(t)
    \right|
    \\
    &\qquad = o\left((V_n^{\omega,\vartheta})^{-1/2}\right).
    \end{aligned}
    \end{equation}
    where $\Phi$ and $\varphi$ are the standard normal distribution function and density, and
    \[
         \kappa_{3,n}^{\omega,\vartheta}
        :=\mathbb E_{\mathbb K_{\vartheta;\omega}}
          \left[
            \left(S_n^{\omega,\vartheta}
              -\mathbb E_{\mathbb K_{\vartheta;\omega}}
                  [S_n^{\omega,\vartheta}]\right)^3
          \right].
    \]
    The same conclusions hold for the sums above viewed as functions of the measurement outcome record under $\mathbb Q _ {\vartheta;\omega}$. 
\end{corollary}
\begin{proof}
    For simplicity, denote $v(\omega) = \operatorname{Var}_{\eta_\omega}(f_{\theta\om})$.
    Measurability of $f$ and $\eta_\omega$ makes $v$ measurable, and boundedness of $f$ gives $v\in L^1(\mathbb P)$.

    We claim that $\mathbb P\{v > 0\} > 0$. 
    Indeed, if $v(\omega) = 0$ then $\mu_\om$ is a point mass, therefore $|\widehat\mu_\om(t)| = 1$ for every $t$. 
    Thus applying, \eqref{eq:qt-nonlattice} at any fixed $t\neq 0$ we get that $\mathbb P\{v > 0 \}> 0 $ and hence $\mathbb E_{\mathbb P}v>0$.     
    
    Now, for each $j\ge 1$ the quenched law of $\rho_j^{\omega,\vartheta}$ dominates $\varepsilon \eta_{\theta^{j-1}\om}$ via the minorization assumption \eqref{eq:qr-one-step}.
    Consequently,
    \begin{align*}
        \operatorname{Var}_{\mathbb K_{\vartheta;\omega}}
        \left(f_{\theta^j\omega}(\rho_j^{\omega,\vartheta})\right)
        &=
        \inf_{a\in\mathbb R} \int \left( f_{\theta^j\omega}(\rho_j^{\omega,\vartheta}) -a\right)^2 \,d \mathbb K_{\vartheta;\omega}\\
        &\ge \varepsilon\inf_{a\in\mathbb R} \int_{\mathbb S_d}  \left( f_{\theta^j\omega}(\rho) -a\right)^2 \,\eta_{\theta^{j-1}\om}(d\rho) \\
        & = \varepsilon\operatorname{Var}_{\eta_{\theta^{j-1}\omega}} (f_{\theta^j\omega})\\
        &= \varepsilon v(\theta^{j-1}\omega).
    \end{align*}
    Therefore using \eqref{eq:variance-comparison} and as $f$ is bounded we obtain the existence of deterministic constants $C_1, C_2>0$, independent of the initial state, such that 
    \begin{equation}
        C_1\varepsilon\sum_{j=0}^{n-1}v(\theta^j\omega)
         \leq V_n^{\omega,\vartheta}
         \leq C_2 n.
    \end{equation}

    Now by ergodicity of $\theta$ and as $v\in L^1(\mathbb P)$, Birkhoff's ergodic theorem gives
    \[
        \frac1n\sum_{j=0}^{n-1}v(\theta^j\omega)
            \longrightarrow\mathbb E_{\mathbb P}v>0
    \]
    for $\mathbb P$-almost every $\omega$. 
    Thus with $\mathbb P$-probability $1$, there is $N(\omega)$ such that the average above is at least $\frac12 \mathbb E_{\mathbb P}v$ for all $n\ge N(\omega)$.
    Thus we obtain that $ c_*n \le V_n^{\omega,\vartheta} \le C_*n $ for sufficiently large $n$ where 
    \[
        c_* = \tfrac12C_1\varepsilon\mathbb E_{\mathbb P}v>0,
         \qquad\text{ and}\qquad C_*=C_2>0.
    \]
    These constants and the threshold $N(\omega)$ are independent of the initial state.

    Under \eqref{eq:qt-dictionary}, the non-lattice hypothesis follows from $(f_\omega)_*m_\omega=\mu_{\theta^{-1}\omega}$ and measure preservation.
    
    Since the variance diverges, Theorem~\ref{LLT rand0}(i) now gives the local limit theorem and the first-order Edgeworth expansion.
    The outcome-space statement follows from \eqref{eq:qt-state-path-law}.
\end{proof}

We emphasize that the conclusion above requires the chosen function $f$ to satisfy the non-lattice condition \eqref{eq:qt-nonlattice}.
While this is a non-trivial assumption, sufficient criteria can be expressed in terms of the posterior states and weights supplied by the instrument process together with the values of $f$ to guarantee this condition. 
For an example, to obtain \eqref{eq:qt-nonlattice}, we use Remark~\ref{rem:nonlattice-criteria} with the additional assumption that $f_{\theta\om}$ is continuous on the compact set $\operatorname{supp}\eta_\omega$ in which case we have that  
\[
     \operatorname{supp}\mu_\omega = f_{\theta\omega}(\operatorname{supp}\eta_\omega).
\]
Thus the criteria in that remark can be checked from the values of $f_{\theta\omega}$ on $\operatorname{supp}\eta_\omega$.

As an example we recall \cite[Proposition~A.2]{traj_paper}, with $L=1$. 
There the outcome space $\Xi$ is finite and there exist $U_1, \ldots, U_q$ disjoint sets of outcomes such that
\[
    \mathfrak I_\omega(\{a\})(X)
      =\operatorname{Tr}(E_{a;\omega}X)\tau_{k;\omega},
\]
for each $a\in U_k$, where $\tau_{k;\omega}\in\mathbb S_d$ and $E_{a;\omega}$ are positive semi-definite matrices, with all data measurable in $\omega$. 
Therefore, upon observing an outcome $a\in U_k$, the posterior state of the system is $\tau_{k;\om}$, independently of the incoming state. 
Furthermore, suppose that 
\[
    A_{k;\omega}:=\sum_{a\in U_k}E_{a;\omega}\geq\alpha_k(\omega)I_d
\]
where $\alpha_k(\omega)\geq0$ are measurable and $\sum_k\alpha_k(\omega)\geq\varepsilon$ for $\mathbb P$-almost every $\omega$.
With this \cite[Proposition~A.2]{traj_paper} gives the one-step minorization given in \eqref{eq:qr-one-step} with 
\[
    \eta_\omega=\sum_kp_k(\omega)\delta_{\tau_{k;\omega}},
    \qquad\text{ and }\qquad 
    p_k(\omega)=\frac{\alpha_k(\omega)}{\sum_\ell\alpha_\ell(\omega)}. 
\]
With this explicit minorization we have that $\mu_\omega$ in \eqref{eq:qt-scalar-law} is exactly 
\[
    \mu_\omega
      =\sum_kp_k(\omega)\delta_{f_{\theta\omega}(\tau_{k;\omega})}.
\]
Set $z_k(\om) := f_{\theta\om}(\tau_{k;\om})$. 
For this mixture, condition~\eqref{eq:qt-nonlattice} is equivalent to requiring that for every $h>0$ we have 
\[
    \mathbb P\left\{
        \om : \exists k, \ell \text{ with } p_k(\om) p_\ell(\om) > 0 \text{ and }
        z_k(\om) - z_\ell(\om) \not\in h \mathbb  Z
    \right\}
    >0
\]
We note that for each fixed $h$ the set above is measurable, being a finite union of measurable events. 
By Remark~\ref{rem:nonlattice-criteria}, a sufficient criterion is that, on an event of positive $\mathbb P$-probability, three values $z_k(\om)$ with positive weights have an irrational ratio of differences. 

However, a more flexible sufficient criterion exists: two values with positive weights can also suffice. 
For example, suppose that for some fixed indices $k,\ell$ the difference $z_k(\om) - z_\ell(\om)$ is nonzero rational on one event of positive probability and irrational on another, with both weights $p_k(\om)$, $p_\ell(\om)$ positive on these events.

\medskip

To obtain an expansion of higher order $r\geq2$ using Theorem~\ref{LLT rand0}(ii), we need an additional condition.
Writing $\widehat\mu_\omega(t)=\int_{\mathbb R}e^{ity}\,\mu_\omega(dy)$, we require, on an event of positive $\mathbb P$-probability,
\[
    |\widehat\mu_\omega(t)|\leq1-c_\omega|t|^{-\ell}
    \qquad(|t|\geq T_\omega),
    \qquad \ell>0,\quad (r-1)\ell<2,
\]
for a fixed integer $r\geq2$ and a deterministic exponent $\ell$.
Here $c_\omega,T_\omega$ are positive and measurable in $\omega$.
Together with the assumptions of the corollary, this yields the expansion of order $r$ under each prescribed initial state for $\mathbb P$-almost every $\omega$.

A convenient sufficient condition is that there exists a measurable event $B\subseteq\Omega$ with $\mathbb P(B)>0$ such that, for every $\omega\in B$, the measure $\mu_\omega$ has a nonzero component absolutely continuous with respect to Lebesgue measure on $\mathbb R$.
Indeed, for each $\om \in B$, write $\mu_\om = \beta_\om \mu_{\om,\mathrm{ac}} + (1 - \beta_\om)\nu_\om$, where $\mu_{\om,\mathrm{ac}}$ is the absolutely continuous part and where $\beta_\om > 0$ is another probability measure. Then, the Riemann-Lebesgue lemma gives $T_\om\ge 1$ such that $|\widehat\mu_{\om,\mathrm{ac}}|\le 1/2$ for $|t|\ge T_\om$.
Therefore, for any $\ell>p$ we have that 
\[
     |\widehat\mu_\omega(t)| \le 1 - \frac{\beta_\om}{2} \le 1 - \frac{\beta_\om}{2}|t|^{-\ell}.
\]

Below we give an example where the apparatus performs the projective measurement with probability $1-\varepsilon$, and with probability $\varepsilon$ it prepares a Haar-distributed pure state $\tau_\psi$ independently of the incoming state.

\begin{example}
\label{eg:qt_1}
    Let $d\ge 2$ and let $\lambda$ be the Haar probability measure on $\mathcal U_d =\{\psi\in\mathbb C^d : \|\psi\| = 1\}$.
    Let $\eta$ be the law of the state $\tau_\psi:=|\psi\rangle\langle\psi|$ under $\lambda$.
    Let $\mathcal A$ be a finite indexed set with $|\mathcal A|\ge 2$, for which we have a mutually orthogonal collection of projections $(P_a)_{a\in\mathcal A}$ such that $\sum_{a\in\mathcal A}P_a = I_d$. 
    Let the outcomes space $\Xi = \mathcal A \sqcup\mathcal U_d$ and let $\varepsilon\in(0,1)$.
    Consider the instrument defined by 
    \[
        \mathfrak I_R(A)(B)
        = (1-\varepsilon)\sum_{a\in A\cap\mathcal A} P_aBP_a
        +\varepsilon \operatorname{Tr}(B)\int_{A\cap\mathcal U_d} \tau_\psi \lambda(d\psi).
    \]
    for $A\subset \Xi$, measurable and $B\in\mathbb M_D(\mathbb C)$.

    We first show that the example satisfies the one-step Doeblin condition in \eqref{eq:qr-one-step}.
    To see this, first let $\mathcal M$ be the posterior state kernel for the projective measurement alone:
    \[
        \mathcal M(\rho,\Gamma)
        = 
        \sum_{a\in\mathcal A, \operatorname{Tr}(P_a\rho)>0} \operatorname{Tr}(P_a\rho) 1_\Gamma\left(\dfrac{P_a\rho P_a}{\operatorname{Tr}{P_a\rho}}\right).
    \]
    Then, we have that 
    \[
        K(\rho, \,\cdot\,) = (1-\varepsilon)\mathcal M(\rho,\,\cdot\,) + \varepsilon\eta,
    \]
    where
    \[
        \eta(\Gamma) = \int_{\mathcal U_d}1_\Gamma(\tau_\psi) \lambda (d\psi).
    \]
    Thus \eqref{eq:qr-one-step} holds.
    Moreover for any self-adjoint $H$ with $H \neq cI_d$, by take $f_H(\rho) = \operatorname{Tr}(H\rho)$ and let $\mu_H := (f_H)_*\eta$.
    If $\{h_1,\ldots , h_d\}$ are the eigenvalues of $H$ then for an orthonormal eigen-basis $\{e_1,\ldots \,e_d\}$ of $H$ with $He_j = h_je_j$ and writing $\psi = \sum_{j=1}^d \psi_je_j$  we get that 
    \[
        f_H(\tau_\psi) = \langle\psi, H\psi\rangle = \sum_{j=1}^d h_j |\psi_j|^2.
    \]
    Since $\psi$ is distributed according to $\lambda$ we have that the vector $(|\psi_1|^2,\ldots,|\psi_d|^2)$ is uniformly distributed on the simplex $\Delta:=\{w =(w_j)_{j=1}^d\in [0,\infty)^d : \sum_{j=1}^d w_j = 1\}$ with respect to the normalized $(d-1)$-dimensional Lebesgue measure; see \cite[Appendix~ eq.A5]{Zyczkowski_2001}. 
    Since $H\neq cI_d$ for any constant $c$ we have that the eigenvalues of $H$ are not all equal; therefore, the affine function $L(H)(w) = \sum_{j=1}^d h_j w_j$ is non-constant on $\Delta$. 
    Thus, by Fubini's theorem, we must have that $\mu_H$ is absolutely continuous: indeed, after relabeling, if necessary, we may assume that $h_1 \neq h_d$ and then for $w\in\Delta$ we have that
    \[
        L(H)(w) = h_d +\sum_{j=1}^{d-1}(h_j - h_d)w_j. 
    \]
    Let $E\subseteq \mathbb R$ be a Lebesgue-null Borel set.
    For each fixed choice of $w_2, \ldots w_{d-1}$, we have that $L(H)$ is an affine map in $w_1$ with a non-zero slope $h_1 - h_d$, yielding that the set of $w_1$ for which $L(H)(w) \in E$ has Lebesgue measure $0$.
    Therefore, Fubini's theorem gives us that $\{w\in\Delta : L(H)(w)\}$ has $0$ $(d-1)$-dimensional Lebesgue measure. Hence $\mu_H(E) = 0$, as well and in particular $\mu_H$ must have infinite support and thus by Remark~\ref{rem:nonlattice-criteria} we have that \eqref{eq:qt-nonlattice} holds for $f_\om = f_H$. 
    Thus Corollary~\ref{cor:llt_qt_1} applies.
    Moreover, by the preceding paragraph before this theorem, we also have that Theorem~\ref{LLT rand0}(ii) holds and gives the Edgeworth expansion of order $r$.
\end{example}


\subsection{Block Minorization and Edgeworth Expansion}

In this section, we apply Theorem~\ref{LLT rand} for disordered quantum trajectories described above. 
For a bounded measurable function $g:\mathbb S_d \to \mathbb C$, write 
\[
    (K_\om g)(\rho) = \int_{\mathbb S_d} g(\sigma) K_\om(\rho, d\sigma)
\]
then for any $L\ge 1$ take 
\[
    K^{(L)}_\om = K_\omega\circ K_{\theta\omega}\circ\cdots\circ K_{\theta^{L-1}\omega}. 
\]
Hence, $K_\omega^{(L)}(\rho,\cdot)$ is the distribution of the posterior state after $L$ measurements, started at $\rho$ in environment $\omega$.
Suppose that there exist a deterministic integer $L\geq1$, a deterministic constant $\varepsilon\in(0,1)$, and a measurable family of probability measures $\eta_\omega^{(L)}$ on $\mathbb S_d$ such that
\begin{equation}
\label{eq:qt-block-minorization}
    K_\omega^{(L)}(\rho,\,\cdot\,)\geq\varepsilon\eta_\omega^{(L)}(\cdot).
\end{equation}
We recall that \cite[Appendix ~A]{traj_paper} gives examples for which such a condition holds with a deterministic length $L$ and is the main assumption used in the results there. 
We take $\Omega_0\in\mathcal F$ to be a set of full $\mathbb P$-probability for which \eqref{eq:qt-block-minorization} holds for all incoming state and take $\Omega_D = \cap_{j\in\mathbb Z}\theta^{-j}\Om_0$. 
Then from \cite[Lemma~4.1]{traj_paper} we have that for every random initial state $\vartheta$ and every $\omega\in\Omega_D$, the posterior chain under $\mathbb K_{\vartheta;\omega}$ consequently satisfies
\begin{equation}
\label{eq:qt-block-mixing}
    \phi_{\vartheta;\omega}(n)
       \leq(1-\varepsilon)^{\lfloor n/L\rfloor}.
\end{equation}
Where  $\phi_{\vartheta;\omega}(n)$ is the $\phi$-mixing coefficient under $\mathbb K_{\vartheta;\om}$ with the forward and backwards filtration generated by the posterior states. 
Therefore, we automatically have the mixing assumption required in Assumption~\ref{A2}. 

To apply Theorem~\ref{LLT rand}, we require the additional conditions from Assumption~\ref{A3}.
Therefore, with the identification that $R_\omega=K_\omega$ and $\cE_\omega=\mathbb S_d$, we impose that:
The space $\Omega$ is topological, $\cF$ is its Borel $\sigma$-algebra, and there is a periodic point $\omega_*\in\Omega$ of period $p\geq1$.
Furthermore, writing $\omega_s^*=\theta^s\omega_*$, we require, for each $s=0,\ldots,p-1$, that
\begin{enumerate}
    \item every open neighborhood $U$ of $\omega_s^*$ has $\mathbb P(U)>0$;
    \item $\theta$ is continuous at $\omega_s^*$;
    \item as $\omega\to\omega_s^*$,
    \[
        \|f_\omega-f_{\omega_s^*}\|_\infty\longrightarrow0,
        \qquad\text{ and }\qquad
        \|K_\omega-K_{\omega_s^*}\|_{\infty\to\infty}\longrightarrow0.
    \]
\end{enumerate}
Here the first norm $\|\,\cdot\,\|_\infty$ is the supremum norm over the state space $\mathbb S_d$ and the operator norm $\|\,\cdot\,\|_{\infty\to\infty}$ is taken on the space of bounded Borel complex-valued functions $B_b(\mathbb S_d)$:
\[
    \|K_\omega-K_{\omega_s^*}\|_{\infty\to\infty}
      =\sup_{\|g\|_\infty\leq1}\sup_{\rho\in\mathbb S_d}
         \left|\int_{\mathbb S_d}g(\sigma)
            \bigl(K_\omega(\rho,d\sigma)-K_{\omega_s^*}(\rho,d\sigma)\bigr)\right|.
\]

For $t\in\mathbb R$ define
\[
     (P_{\omega,t}g)(\rho)
       =\int_{\mathbb S_d}
          e^{itf_{\theta\omega}(\sigma)}g(\sigma)\,K_\omega(\rho,d\sigma),
\]
and take
\[
    Q_t=P_{\omega_*,t}\circ\cdots\circ P_{\theta^{p-1}\omega_*,t}.
\]
With these notations, a further condition that is required to apply Theorem~\ref{LLT rand} is that 
\begin{equation}\label{eq:qt-periodic-spectral}
    r(Q_t)<1\qquad\text{for every }t\in\mathbb R\setminus\{0\},
\end{equation}
where $r(Q_t)$ denotes the spectral radius on $B_b(\mathbb S_d)$.

\medskip 

Let $f:\Omega\times\mathbb S_d\to\mathbb R$ be a jointly measurable and uniformly bounded function, and let $S_n^{\omega,\vartheta}$ and $V_n^{\omega,\vartheta}$ be as in \eqref{eq:qt-sum}. 
Under \eqref{eq:qt-block-minorization}, \cite[Theorem~3]{traj_paper} gives a deterministic constant $\Sigma_f^2\geq0$, independent of the initial state, such that, for every measurable initial state $\vartheta$,
\begin{equation}
\label{eq:qt-variance-limit}
    \frac{V_n^{\omega,\vartheta}}{n}\longrightarrow\Sigma_f^2
    \qquad\text{for $\mathbb P$-almost every $\omega$}.
\end{equation}
This limit also holds for the arrival-state sums in \eqref{eq:qt-sum}, since changing the two bounded endpoint terms does not change the variance divided by $n$ under the bound $V_n=O(n)$.
The constant $\Sigma_f^2$ may be zero; its positivity is required below.

\begin{corollary}
\label{cor:qt-block-local}
    Suppose \eqref{eq:qt-block-minorization}, the periodic-support and continuity conditions above, and assume that \eqref{eq:qt-periodic-spectral} holds.
    Assume also that $\Sigma^2_f > 0$.
    Then, for any random initial state $\vartheta:\Omega\to\mathbb S_d$ and $\mathbb P$-almost every $\om$, the sums $S_n^{\omega,\vartheta}$ satisfy the non-lattice local central limit theorem of Definition~\ref{dfn:LCLT} and the first-order Edgeworth expansion \eqref{eq:qt-first-edgeworth} under $\mathbb K_{\vartheta;\omega}$.
\end{corollary}
\begin{proof}
    By \eqref{eq:qt-block-mixing}, the posterior chain under  $\mathbb K_{\vartheta;\omega}$ satisfies the uniform exponential $\phi$-mixing bound for $\mathbb P$-almost every $\omega$, in Assumption~\ref{A2}.
    The remaining assumptions give Assumption~\ref{A3}, the periodic spectral condition, and positive asymptotic variance $\sigma^2 =\Sigma^2_f>0$. 
    The forward estimates of Theorem~\ref{LLT rand}, applied with $R_\omega=K_\omega$ and initial state  $\vartheta(\omega)$, therefore yield both conclusions.
\end{proof}

We note that the environmental conditions required in the result above can be checked using Examples~\ref{ex:periodic-subshift} and~\ref{ex:periodic-markov}.
To verify \eqref{eq:qt-block-minorization} we may use \cite[Proposition~A.2]{traj_paper} or other criterion from \cite[Appendix~A.]{traj_paper}. 
To verify \eqref{eq:qt-periodic-spectral}, we may use Remark~\ref{rem:periodic-spectral-criterion}. 

\medskip

Now we give an example.
In the example below, an external classical controller chooses the instrument at each time independently of the measurement outcome, and the measurement record includes the instrument type used as well as the recorded outcome.  

\begin{example}
    Denote by $\mathfrak I_R$ the instrument from Example~\ref{eg:qt_1} and denote by $\Xi_R$ the outcome space there.
    Denote by $K_R$ the kernel from  Example~\ref{eg:qt_1} 
    Now for $m\ge 1$ let $\mathfrak I_{A_1}, \ldots, \mathfrak I_{A_m}$ be some finite collection of instruments with the associated Borel outcome spaces $\Xi_{A_j}$ for each $j\in\{1,\ldots,m\}$.
    Denote by $K_{a_j}$ the posterior kernel for the instrument $\mathfrak I_j$.
    We now treat the collection of instruments $\{\mathfrak I_R,\mathfrak I_{A_1}, \ldots, \mathfrak I_{A_m}\}$ on
    \[
        \Xi = \Xi_R\sqcup \Xi_{A_1}\sqcup\cdots\sqcup\Xi_{A_m}
    \]
    by extending each instrument trivially outside of its own component. 

    Fix $L\ge 2$ and $q\in(0,1)$.
    Let $(C_n)_{n\ge 0}$ be a Markov process on some auxiliary probability space $(\widetilde\Omega,\widetilde{\mathcal F},\widetilde{\mathbb P})$ with the state space $\{0,\ldots,L-1\}$ and the transition probabilities given by 
    \[
        p(k,0)=q,\qquad p(k,k+1)=1-q\quad \text{for} \quad 0\le k<L-1, \qquad \text{and} \qquad p(L-1,0)=1.
    \]
    This Markov chain is irreducible and aperiodic with the stationary distribution 
    \[
        \pi_k = \dfrac{q(1-q)^k}{1 - (1-q)^L}\qquad 
        \text{ for }
        \qquad
        k\in\{0,\ldots ,L-1\}.
    \]
    At time $n$ we select the instrument $Z_n = R$ is $C_{n+1} = 0$; else we select an instrument uniformly from $\{A_1,\ldots,A_m\}$. We take $C_0\sim\pi$. 
    Now let 
    \[
        \Om = \{\om\in\{R, A_1, \ldots ,A_m\}^{\mathbb Z} : \text{ every block of length $L$ contains $R$ }\},
    \]
    equipped with the Borel $\sigma$ algebra under the product topology. 
    For integers $m\le n$ and symbols $z_m,\ldots , z_n\in \{R, A_1, \ldots, A_m\}$ we prescribe the probabilities 
    \[
        \mathbb P \{ \omega_a=z_a,\ldots,\omega_b=z_b\} = \widetilde{\mathbb P}\{Z_0=z_a,\ldots,Z_{b-a}=z_b\}.
    \]
    Stationarity makes these prescriptions consistent, so the Kolmogorov extension theorem gives a probability measure on $\Om$ which we denote by $\mathbb P$.
    Taking the left shift by $(\theta\omega)_j=\omega_{j+1}$, we have that $(\Om, \mathcal F, \mathbb P,\theta)$ is ergodic. 
    This space $\Om$ is a topologically mixing subshift of finite type: the words made from symbols $\{R, A_1, \ldots, A_m\}$ that are not allowed are the words of length $L$ that do not contain the symbol $R$, and any allowed words can be joined by an arbitrarily long string of the symbol $R$. 
    The schedule is a stationary Markov process of order at most $L-1$, with positive conditional probabilities for every admissible next symbol. 
    Consequently, its law $\mathbb P$ is a Gibbs measure for a locally constant potential $\log p(\omega_{L-1}\mid\omega_0,\ldots,\omega_{L-2})$.
    
    We prove that this instrument process satisfies \eqref{eq:qt-block-minorization}.
    Set $\mathfrak I_\omega=\mathfrak I_{\omega_0}$ and $K_\omega=K_{\omega_0}$.
    Now if $s(\om)$ is the last occurrence of $R$ in $\om_0,\ldots, \om_{L-1}$ then we have that 
    \[
        K_\omega^{(L)}(\rho,\cdot)\ge\varepsilon\eta_\omega^{(L)},
    \]
    where
    \[
        \eta_\omega^{(L)}:=\eta K_{\omega_{s+1}}\cdots K_{\omega_{L-1}}. 
    \]
    Thus \eqref{eq:qt-block-minorization} holds and \eqref{eq:qt-block-mixing} supplies Assumption~\ref{A2}.

    Now, for a self-adjoint matrix $H$ with $H\neq c_d$ for any constant $c$, take $f_\omega = f_H$ and $\mu_H$  as in the previous example. 
    We note that Assumption~\ref{A3} holds for this space with the fixed point $\omega_* = (\ldots , R, R, R,\ldots ) = R^{\mathbb Z} \in \Om$, since $K_\om$ depends only on $\om_0$ (and as $f$ is constant). 
    Furthermore, at $ R^{\mathbb Z}$ we have that $Q_tg = K_R(e^{itf_Hg})$ and thus by Remark~\ref{rem:periodic-spectral-criterion} gives 
    \[
         \|Q_t^2\|_{\infty\to\infty} \le 1-\varepsilon^2(1 - |\widehat\mu_H(t)|) < 1
    \]
    for $t\neq 0$. 

    Now we verify that $\Sigma^2_H > 0$.
    For $n\ge 2$, let
    \[
        M_n(\om) := \sum_{j=0}^{n-2}
        \mathbf 1_{\om_j = \om_{j_+1} = R}.
    \]
    Since $\mathbb P\{\om_0 = \om_1 = R\}>0$ we have from the Birkhoff's ergodic theorem that
    \[
        M_n(\om) \to \mathbb P\{\om_0 = \om_1 = R\}: = p_{RR} >0
    \]
    for $\mathbb P$-almost every $\om$. 
    Now denote by $N_n(\om)$ be the number of disjoint $RR$ \emph{pairs} selected from $\omega_0,\ldots,\omega_{n-1}$ as follows: If $j\le n-2$ and $\om_j - \om_{j+1} = R$ we select $(j,j+1)$ and continue from $j+2$, else we advance to $j+1$ and stop counting pairs when we reach $j\ge n-1$. 
    Here no position belongs to two selected pairs. 
    This gives us that $M_n(\om) \le 2N_n(\om)$, and thus 
    \[
        \liminf_{n\to\infty}\frac{N_n(\om)}{n} > 0
    \]
    for $\mathbb P$-almost every $\om$. 
    Now for any random initial state $\vartheta$, the Markov property gives that 
    \[
         \mathbb E_{\mathbb K_{\vartheta;\omega}} [e^{itS_n^{\omega,\vartheta}}]
         =
         \left(
            P_{\omega,t}\circ\cdots\circ P_{\theta^{n-1}\omega,t}\mathbf 1
         \right) (\vartheta(\omega))
    \]
    where $P_{\omega,t}g=K_\omega(e^{itf_H}g)$. 
    Therefore, we have that 
    \[
         \left|
            \mathbb E_{\mathbb K_{\vartheta;\omega}}
             [e^{itS_n^{\omega,\vartheta}}]
        \right|
        \le \left(1-\varepsilon^2(1 - |\widehat\mu_H(t)|)\right)^{N_n(\om)}.
    \]
    Now,  for any random variable $Y\in L^3$, 
    \[
        |\mathbb E[e^{itY}]| 
        =
        1-\tfrac12\operatorname{Var}(Y)t^2+o(t^2) 
    \]
    as $t\to0$. Set $v_H=\operatorname{Var}_\eta(f_H)>0$.
    Therefore $V_n^{\omega,\vartheta}\ge\varepsilon^2v_HN_n(\omega)$ and thus dividing by $n$ and applying \eqref{eq:qt-variance-limit} we get that 
    \[
        \Sigma_H^2
        \ge\varepsilon^2v_H
            \liminf_{n\to\infty}\frac{N_n(\omega)}n
        \ge\frac{\varepsilon^2v_Hp_{RR}}2>0.
    \]
    
    Corollary~\ref{cor:qt-block-local} now gives the non-lattice local central limit theorem and the first-order Edgeworth expansion.
    Moreover, absolute continuity of $\mu_H$ and the Riemann--Lebesgue lemma give $a(t)\le1-\varepsilon^2/2$ for all sufficiently large $|t|$.
    For every fixed $r\ge2$, Theorem~\ref{Edge2} therefore applies with $m_0=2$ and any $0<\ell<2/(r-1)$, giving an Edgeworth expansion of order $r$.
    These conclusions hold for every measurable initial state and $\mathbb P$-almost every environment.
\end{example}


\bibliography{ref}

\end{document}